\documentclass[]{article}
\usepackage{geometry}
\usepackage{amsmath,amssymb,amsfonts,amsthm,amsopn,dsfont, mathrsfs}
\usepackage{bm}
\usepackage{graphicx}
\usepackage{float}
\usepackage[backend=bibtex,maxnames=5]{biblatex}
\newtheorem{prop}{Proposition}
\newtheorem{lemma}{Lemma}
\newtheorem{rem}{Remark}
\newtheorem{assum}{Assumption}

\usepackage{standalone}
\usepackage{tikz}
\usepackage{pgfplots}
\usepackage{pgfplotstable}
\usepackage[title]{appendix}
\pgfplotsset{compat=newest}
\usepackage{subcaption}
\usepackage{colortbl}
\usepackage{booktabs}
\usepackage{stmaryrd}

\usepackage[T1]{fontenc}

\usepackage{xcolor}
\newcommand{\OmegaT}{\Omega_{_\mathcal{T}}}
\newcommand{\flux}{\bm{\sigma}_{h}}

\newcommand{\fluxa}{\bm{\sigma}_h^{\bm{a}}}

\newcommand{\tr}[2]{{#1}|_{{#2}}}
\newcommand{\T}[1]{{\mathcal{T}_{#1}}}
\newcommand{\A}[1]{{\mathcal{A}_{#1}}}

\newcommand{\N}[1]{{\mathscr{N}_#1}}
\newcommand{\Hdiv}[1]{\bm{H}(\text{div},#1)}
\newcommand{\estim}[1]{\mathcal{E}_{#1}}
\newcommand{\Mha}{\bm{M}_h^{\bm{a}}}
\newcommand{\Qha}{Q_h^{\bm{a}}}
\newcommand{\wa}{{\omega_{\bm{a}}}}
\newcommand{\was}{{\omega}_{\bm{a}}^\star}
\newcommand{\psia}{{\psi_{\bm{a}}}}
\newcommand{\lam}{\lambda_h^{\bm{a}}}

\newcommand{\gr}{\tilde{\gamma}}
\newcommand{\gs}{{\gamma^\star}}

\newcommand{\Fs}{\mathcal{F}^\star}
\newcommand{\Fr}{\tilde{\mathcal{F}}}
\newcommand{\us}{u^\star}
\newcommand{\uhs}{u_h^\star}
\newcommand{\partwa}[1]{{\partial \omega_{\bm{a}}^{#1}}}

\newcommand{\Omegas}{{\Omega^\star}}

\newcommand{\polydeg}{{r}}

\begin{document}
	%opening
	\title{Adaptive refinement in defeaturing problems via an equilibrated flux \textit{a posteriori} error estimator}
	\author{Annalisa Buffa\thanks{MNS, Institute of Mathematics, École Polytechnique Fédérale de Lausanne, Switzerland \& IMATI CNR - Via Ferrata 5, 27100 Pavia  ({annalisa.buffa@epfl.ch})} \and Denise Grappein\thanks{MOX, Department of Mathematics, Politecnico di Milano, via Bonardi 9, 20133, Member of GNCS INdAM Group ({denise.grappein@polimi.it})} \and Rafael V\'azquez\thanks{Departamento de Matem\'atica Aplicada, Universidade de Santiago de Compostela, Spain \& Galician Centre for Mathematical Research and Technology (CITMAga), Santiago de Compostela, Spain ({rafael.vazquez@usc.es})}}

	\maketitle
	\begin{abstract}
		An adaptive refinement strategy, based on an equilibrated flux \textit{a posteriori} error estimator,
	is proposed in the context of defeaturing problems. Defeaturing consists of removing features from complex domains to simplify mesh generation and reduce the computational cost of simulations. It is a common procedure, for example, in computer aided design for simulation-based manufacturing. However, depending on the problem at hand, geometrical simplification may significantly deteriorate the accuracy of the solution. The proposed adaptive strategy is hence twofold: starting from a defeatured geometry, it iteratively performs both standard mesh refinement and geometrical refinement by selecting, at each step, which features must be reintroduced to significantly improve accuracy. Similar adaptive strategies have been previously developed using residual-based error estimators within an IGA framework. Here, instead, we extend a previously developed equilibrated flux \textit{a posteriori} error analysis, designed for standard finite element discretizations, to make it fully applicable within the adaptive procedure. In particular, we address the assembly of the equilibrated flux estimator in presence of elements trimmed by the boundary of included features, adopting an unfitted finite element strategy to handle feature inclusion. The resulting estimator allows us to bound both the defeaturing and the numerical sources of error, with additional contributions accounting for the weak imposition of boundary conditions.
\end{abstract}

	\section{Introduction}
Defeaturing consists of the simplification
of a geometry by removing features that are considered not relevant for the approximation
of the solution of a given PDE. It is a fundamental process to reduce the computational effort when repeated simulations on complex geometries are required, in particular when the computational domain is characterized by the presence of features of different scales and shapes. This is the case, for example, of simulation-based manufacturing. However, identifying which features actually have a negligible impact on the solution may be not trivial. Historically, the problem has been approached by exploiting some \emph{a priori} knowledge of the domain, of the materials and of the problem at hand (\cite{FINE00,FOUCAULT2004,THAKUR2009}). Nevertheless, an \emph{a posteriori} estimator becomes necessary when the industrial design process needs to be automatized, and several examples can be found in the literature. In \cite{Ferrandes2009} the defeaturing error is assumed to be concentrated on the feature boundary and feature-local problems are solved in order to estimate it. Alternative approaches are instead proposed in \cite{CHOI2005,GOP2007,GOP2009,SOK1999,TUR2009}, resorting to the concept of feature sensitivity analysis. In \cite{LI2011b,LI2011,LI2013,LI2013b} an \emph{a posteriori} estimator is built by reformulating the defeaturing error as a modeling error, whereas \cite{TANG2013} resorts to the reciprocal theorem, stating flux conservation in the features.

{Recently, a precise framework for analysis-aware defeaturing for the Poisson equation was introduced  in \cite{BCV2022}, proposing an estimator which explicitly depends on the feature size. A geometric adaptive strategy based on this estimator was devised in \cite{AC_2023}, enabling, at each step, the selective inclusion of features according to their impact on the solution accuracy. In \cite{BCV2022_arxiv}, geometric adaptivity was integrated with standard mesh adaptivity within an IGA discretization framework, by coupling the defeaturing estimator of \cite{BCV2022} with a residual-based \textit{a posteriori} estimator for the numerical error. These works were the starting point for the investigation carried out in
\cite{BCGVV}, where the residual-based estimator was replaced by an equilibrated flux estimator, bounding the numerical component of the error with no unknown constant  (\cite{AINSWORTH1997,DEST_MET1999,BraessScho2008,Luce_Wohlmuth}). Employing the equilibrated flux also within the defeaturing component made the estimator particularly suitable for standard finite element discretizations, in which the numerical flux is typically discontinuous across element faces. In \cite{BCGVV} the equilibrated flux is reconstructed on a mesh that is completely unaware of the presence of the features, using the local equilibration procedure proposed in \cite{BraessScho2008, EV2015}, and integrals on the feature boundary, which are required to evaluate the defeaturing component of the estimator, can always be computed, regardless of the intersections with the mesh elements.

The analysis presented in \cite{BCGVV} was intended as a preliminary step towards the use of the equilibrated flux estimator in an adaptive strategy similar to the one devised in \cite{BCV2022_arxiv}, combining mesh and geometric refinement. However, the \textit{a posteriori} error analysis carried out in \cite{BCGVV} remains valid only until the first feature is marked for inclusion. Once a feature is included, the mesh becomes non-conforming to an actual domain boundary, and the flux reconstruction and the estimator assembly need to be addressed in presence of trimmed mesh elements. Alternatively, repeated remeshing of the computational domain would be required after each feature inclusion, which is, of course, not a viable option from a computational cost perspective.

In this work we extend the analysis of \cite{BCGVV} to build an \textit{a posteriori error estimator} on partially defeatured geometries, where selected features have already been reintroduced into the domain. We focus on the negative Neumann feature case, namely small holes on which Neumann boundary conditions are imposed. Observing that the Neumann boundary condition is essential for the reconstructed flux, we adapt the {unfitted finite element} strategy proposed in \cite{puppi2021cut} to weakly enforce the boundary condition on the immersed feature boundary during flux reconstruction. The resulting weak equilibration introduces additional contributions to the estimator. Nevertheless, the reconstruction is designed so that these additional terms vanish on elements that are not directly intersected by the feature boundary. The final \emph{a posteriori} error estimator retains the decomposition into a defeaturing and a numerical component: the former bounds the error due to the features that are not yet included into the geometry, while the latter consists of the standard constant-free equilibrated flux estimator on uncut elements, together with additional contributions on trimmed elements. These contributions account for the weak imposition of Neumann boundary conditions and the resulting loss of local mass conservation. {The analysis focuses on the reliability of the proposed estimator. Efficiency is not addressed here, as the error and the estimator are defined on domains with a different number of features, so that a standard efficiency-proof setting does not apply. Numerical experiments show, however, that the error is not significantly overestimated and that the estimator effectively drives the adaptive strategy.} }

The manuscript is organized as follows: in Section \ref{sec:not_mod_prob} we introduce notation and the model problem, while details about the flux reconstruction procedure are reported in Section \ref{sec:flux}. Section \ref{sec:aposteriori} is devoted to the derivation and analysis of an \emph{a posteriori} error estimator for the overall error and for the proof of its reliability, while the adaptive procedure based on the estimator is detailed in Section \ref{sec:adaptive}. Finally, in Section \ref{sec:num_exp} some numerical tests are proposed, in order to validate the proposed estimator.
	\section{Notation and model problem}\label{sec:not_mod_prob}
	In the following, given any open $k$-dimensional manifold $\omega\subset \mathbb{R}^d$, $d=2,3$ and $k\leq d$, we denote by $|\omega|$ the measure of $\omega$, by $(\cdot,\cdot)_\omega$ the $L^2$-inner product on $\omega$ and by $||\cdot||_\omega$ the corresponding norm. If $k<d$, then $\langle \cdot,\cdot \rangle_\omega$ stands for a duality paring on $\omega$. We denote by $\partial \omega$ the boundary of $\omega$ and, given $\varphi\subseteq \partial \omega$ and $z\in H^{\frac{1}{2}}(\varphi)$, we define
	$$H_{z,\varphi}(\omega)=\lbrace v \in H^1(\omega):~\tr{v}{\varphi}=z \rbrace.
	$$
	For future use we define the quantity
	\begin{equation}
	c_\omega:=\begin{cases}
	\max(-\log(|\omega|),\zeta)^{\frac{1}{2}}&\text{if } k=1,~ d=2\\
	1 &\text{if } k=2, ~d=3,
	\end{cases}\label{c_omega}
	\end{equation}
	where $\zeta \in \mathbb{R}$ is the unique solution of $\zeta=-\log(\zeta)$. 
	Finally, we use the symbol $\lesssim$ to denote any inequality which does not depend on the size of the considered domains, but which can depend on their shape.
	 
\begin{figure}
	 	\centering
	 	\includegraphics[width=0.2\linewidth]{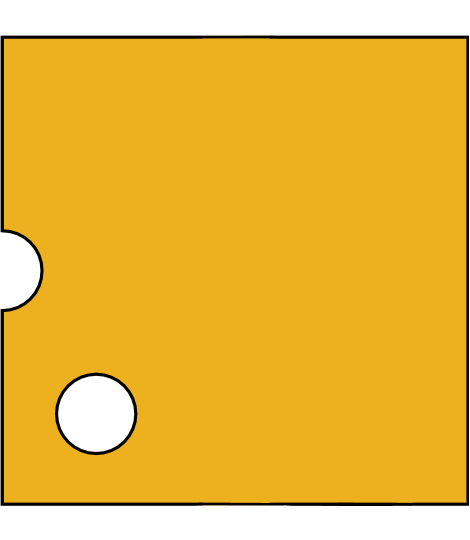}\hspace{4cm}%
	 	\includegraphics[width=0.2\linewidth]{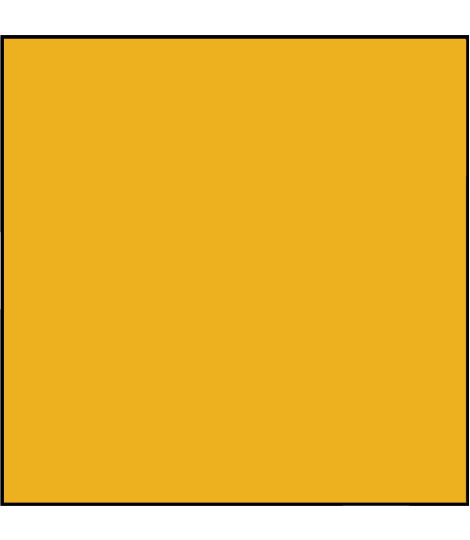}
	 	\caption{Left, domain with two negative features; right, corresponding defeatured geometry.}
	 	\label{fig:posandneg}
	 \end{figure}

	Let $\Omega \subset \mathbb{R}^d$ be an open domain with boundary $\partial \Omega$, characterized by a finite set of open Neumann negative features $\mathcal{F}$. We say that $F \in \mathcal{F}$ is a negative feature of $\Omega$ if it is a geometric detail of smaller scale such that $(\overline{F}\cap \overline{\Omega})\subset \partial \Omega$. In practice, a negative feature corresponds to a portion of the domain in which some material has been removed, internally or at the boundary (see Figure~\ref{fig:posandneg}-left). We call $F\in \mathcal{F}$ a Neumann feature if a Neumann boundary condition is imposed on $\partial F \cap \partial \Omega$. We refer to Remark \ref{rem:pos} to some comments on different types of features and boundary conditions. 
	For the sake of simplicity in the analysis that follows, we assume that both $\Omega$ and each $F\in \mathcal{F}$ are Lipschitz domains, and we denote by $\partial F$ the boundary of a generic feature $F$.	{We assume that the features are pairwise disjoint, i.e. that for every $F,F'\in \mathcal{F}$, $F\neq F'$, $\overline{F}\cap \overline{F'}=\emptyset$. This assumption is not particularly restrictive, as if two features intersect each other we can simply redefine them as a single feature.}
	 
	 The simplified geometry obtained by neglecting all the features is called \emph{defeatured geometry}. We denote it by $\Omega_0$ and, in practice, we obtain it by filling with material the negative features (see Figure \ref{fig:posandneg}-right), i.e. \begin{equation}
	 \Omega_0:=
	 \text{int}\big(\overline{\Omega}\cup\lbrace\overline F,~F \in \mathcal{F}\rbrace\big)\label{Omega0}.
	 \end{equation} For simplicity, we assume $\Omega_0$ to be a Lipschitz domain as well. 

	 Let $\partial \Omega=\Gamma_{\mathrm{D}}\cup\Gamma_{\mathrm{N}}$, with $\Gamma_{\mathrm{D}}\cap \Gamma_{\mathrm{N}}=\emptyset$, and, as in \cite{BCGVV}, we assume that $\partial F\cap \Gamma_D=\emptyset$, $\forall F \in \mathcal{F}$, which means that all the features lie on a Neumann boundary. Finally, let $$\gamma_0:=\bigcup_{F \in  \mathcal{F}}(\partial F \setminus \overline{\Gamma_{\mathrm{N}}})\subset \partial \Omega_0, \quad \text{and} \quad \gamma:=\bigcup_{F \in \mathcal{F}}(\partial F\setminus \overline{\gamma_0})\subset \partial \Omega.$$ We set $\gamma_{_F}=\tr{\gamma}{F}$, such that $\gamma=\bigcup_{F \in  \mathcal{F}}\gamma_{_F}$.\\
	 \begin{figure}
	 	\centering
	 	\includegraphics[width=0.465\linewidth]{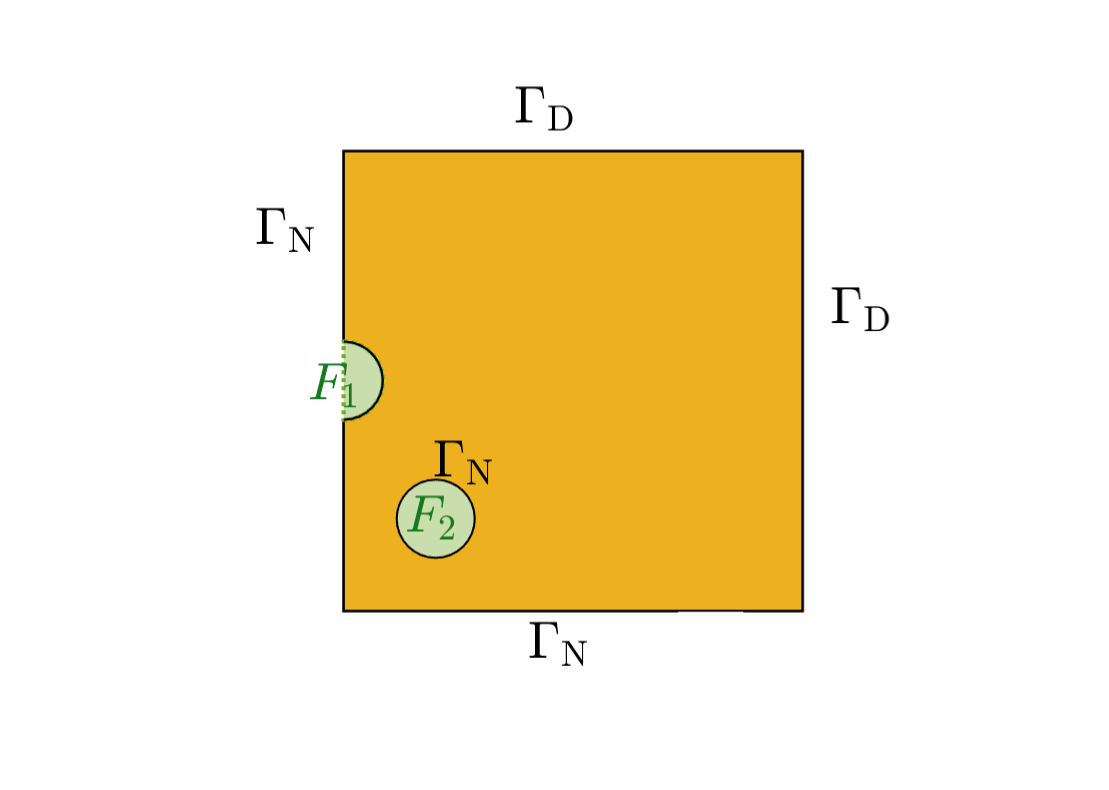}%
	 	\includegraphics[width=0.45\linewidth]{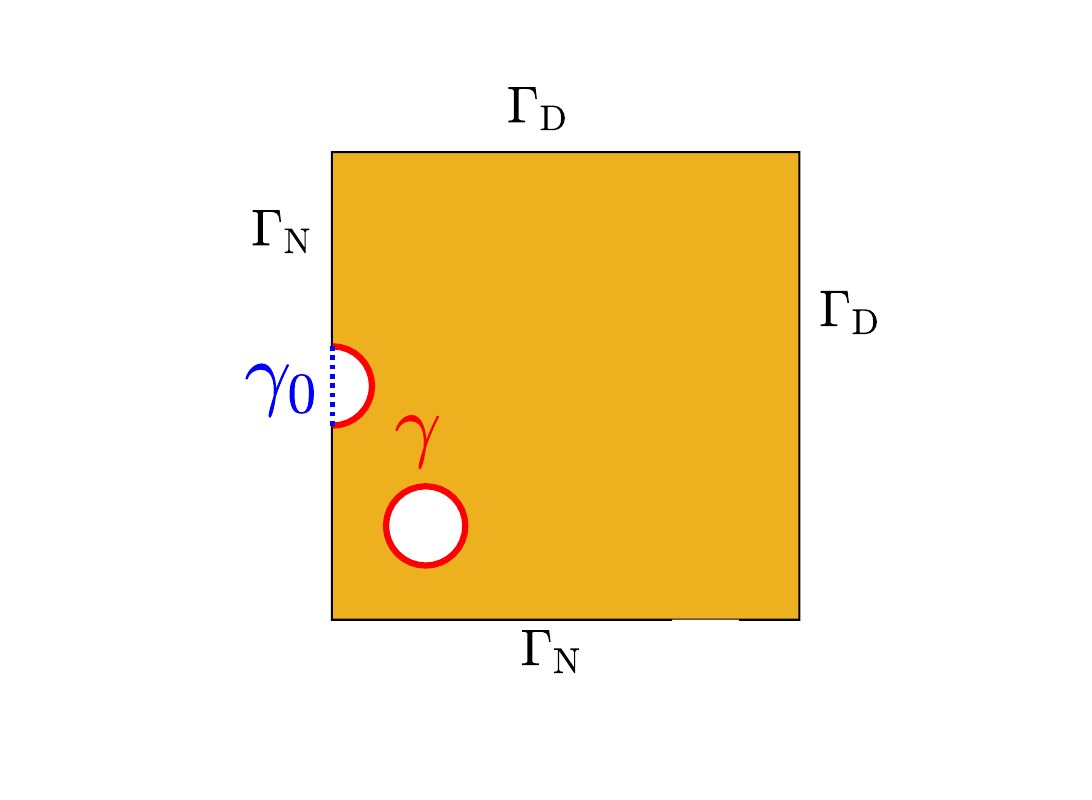}\vspace{-0.8cm}
	 	\caption{Domain with two negative features and corresponding boundary nomenclature.}
	 	\label{fig:gamma_gamma0}
	 \end{figure}
	 
	{We now do two further assumptions on the geometry of the features, similarly to \cite{BCV2022}:}
	{\begin{assum}[Non-accumulation]\label{assum1}
		For every $F \in \mathcal{F}$ there exists a sub-domain $\Omega^F\subset \Omega$ such that $\gamma_{_F}\subset \partial \Omega^F$, and
			the constants $$C_{\mathrm{vol}}:=\max_{F \in \mathcal{F}} \frac{|\Omega|}{|\Omega^F|} \quad \text{and} \quad N_\mathrm{ov}:=\max_{\bm{x} \in \overline{\Omega}}\sum_{F \in \mathcal{F}}\chi_{\overline{\Omega^F}}(\bm{x})$$ are bounded independently of the size of features, where $\chi_{\overline{\Omega^F}}$ denotes the indicator function of $\overline{\Omega^F}$.
	\end{assum}}
{\begin{assum}[Extension operator]\label{assum2}
		There exists an extension operator $\mathbf{E}^{0}:H_{0,\Gamma_D}^1(\Omega)\rightarrow H_{0,\Gamma_D}^1(\Omega_0)$ such that 
		\begin{equation}
			||\nabla \mathbf{E}^{0}(v)||_{\Omega_0}\leq C_{\mathrm{stein}}||\nabla v||_{\Omega}\quad \forall v \in H_{0,\Gamma_D}^1(\Omega) \label{stein}
		\end{equation}
		where the positive constant $C_{\mathrm{stein}}$ is independent of the size and of the number of features.
\end{assum}}
{Let us remark how Assumption~\ref{assum1} excludes the possibility of an increasing number of arbitrarily small features accumulating in the same region of $\Omega$, such as in fractal geometries. Indeed it implies that the number of features is bounded by the quantity $C_{\rm{vol}}N_{\rm{ov}}$. The subdomains $\Omega^F$ are only used to count features and have no relation to any neighborhood used to construct $\mathbf{E}^0$, which may instead be of size comparable to that of the features. Concerning Assumption~\ref{assum2}, an extension operator satisfying \eqref{stein} is built in \cite{Sauter} for a large class of domains, and in particular for domains characterized by separated geometrical details of smaller scale, for which it is required that the mutual distance of the features is at least proportional to their size. The classical construction of \cite{Stein1970} applies instead to any Lipschitz domain, but its constant depends in general on the size of the features. }
	
	Let us choose as a model problem the Poisson problem on $\Omega$:
	\begin{equation}
	\begin{cases}
	-\Delta u =f & \text{in}~ \Omega\\
	u=g_\mathrm{D} &\text{on }  \Gamma_{\mathrm{D}}\\
	\nabla u \cdot \bm{n}=g &\text{on }  \Gamma_{\mathrm{N}},
	\end{cases} \label{eq:strong_omega}
	\end{equation}
	with $\bm{n}$ being the unitary outward normal of $\Omega$. {For the sake of simplicity we here assume $\Gamma_{\mathrm{D}}\neq \emptyset$. Data compatibility conditions would be required in case of a pure Neumann problem.} 
	The variational formulation of Problem \eqref{eq:strong_omega} reads: \emph{find $u \in H^1_{g_\mathrm{D},\Gamma_{\mathrm{D}}}(\Omega)$ that satisfies, $\forall v \in H^1_{0,\Gamma_{\mathrm{D}}}(\Omega)$}
	\begin{equation}
	(\nabla u,\nabla v)_{\Omega}=(f,v)_{\Omega}+\langle g,v\rangle_{\Gamma_{\mathrm{N}}}. \label{eq:weak_omega}
	\end{equation}
	
On the defeatured geometry $\Omega_0$ we consider instead the \emph{defeatured problem} 
	\begin{equation}
	\begin{cases}
	-\Delta u_0 =f & \text{in}~ \Omega_0\\
	u_0=g_\mathrm{D} &\text{on }  \Gamma_{\mathrm{D}}\\
	\nabla u_0\cdot \bm{n}_0=g &\text{on }  \Gamma_{\mathrm{N}}\setminus \gamma\\
		\nabla u_0\cdot \bm{n}_0=g_0 &\text{on }  \gamma_0,
	\end{cases} \label{eq:strong_omega0}
	\end{equation}
	where $\bm{n}_0$ is the unitary outward normal of $\Omega_0$ and, by an abuse of notation, $f \in L^2(\Omega_0)$ is a suitable $L^2$-extension of $f \in L^2(\Omega)$.
	The variational formulation of Problem \eqref{eq:strong_omega0} reads: \emph{find $u_0 \in H^1_{g_\mathrm{D},\Gamma_{\mathrm{D}}}(\Omega_0)$ that satisfies, $\forall v \in H^1_{0,\Gamma_{\mathrm{D}}}(\Omega_0)$}
	\begin{equation}
	(\nabla u_0,\nabla v)_{\Omega_0}=(f,v)_{\Omega_0}+\langle g,v\rangle_{\Gamma_{\mathrm{N}}\setminus \gamma}+\langle g_0,v\rangle_{\gamma_0}. \label{eq:weak_omega0}
	\end{equation}
	
		Let $\mathcal{T}_h^0$ be a non-degenerate simplicial mesh covering $\Omega_0$, consisting of open simplices $K$, such that $\overline{\Omega_0}=\bigcup_{K \in \mathcal{T}_h^0}\overline{K}$. Hereby, we suppose that the mesh faces match with the
	boundaries $\Gamma_\mathrm{D}$, $\Gamma_{\mathrm{N}}^0:=(\Gamma_{\mathrm{N}}\setminus \gamma) \cup\gamma_0$, but we are not asking $\mathcal{T}_h^0$ to be conforming to the boundaries in $\gamma$.
	Let us introduce the space
	\begin{equation}Q_h(\Omega_0)=\mathcal{P}_\polydeg(\mathcal{T}_h^0):=\left\lbrace q_h\in L^2(\Omega_0):~\tr{q_h}{K}\in \mathcal{P}_\polydeg(K),~\forall K \in \mathcal{T}_h^0 \right\rbrace,\label{eq:Qh0}
	\end{equation} 
	with $\mathcal{P}_\polydeg(K)$ denoting the space of polynomials of degree at most $\polydeg$ on $K \in \mathcal{T}_h^0$, with $r\geq 1$, and let
	$$V_h^0(\Omega_0)=\lbrace v \in \mathcal{C}^0(\overline{\Omega_0})\cap Q_h(\Omega_0): \tr{v}{\Gamma_{\mathrm{D}}} =0\rbrace$$
	$$
	V_h(\Omega_0)=\lbrace v \in \mathcal{C}^0(\overline{\Omega_0})\cap Q_h(\Omega_0): \tr{v}{\Gamma_{\mathrm{D}}} =g_\mathrm{D}\rbrace.$$
	For the sake of simplicity, we assume $f \in Q_h(\Omega_0)$.
	Similarly, we assume
	\begin{equation}
	g_\mathrm{N}:=\begin{cases} g &\text{on } \Gamma_{\mathrm{N}}\setminus \gamma\\
	g _0&\text{on } \gamma_0\end{cases} \label{eq:gN}
	\end{equation}
	and $g_D$
	to be piecewise linear polynomials on the partition induced by $\mathcal{T}_h^0$ on $\Gamma_\mathrm{N}^0$ and $\Gamma_{\mathrm{D}}$, respectively. The discrete version of the defeatured problem can then be written as:
	\emph{find $u_h^0 \in V_h(\Omega_0)$ such that}
	\begin{equation}
	(\nabla u_h^0,\nabla v_h)_{\Omega_0}=(f,v_h)_{\Omega_0}+\langle g_\mathrm{N},v_h\rangle_{\Gamma_{\mathrm{N}}^0}\quad \forall v_h \in V_h^0(\Omega_0). \label{eq:prob_num_def}
	\end{equation}

\begin{figure}
		\centering
		\begin{subfigure}{.3\textwidth}
			\centering
			\includegraphics[width=0.75\linewidth]{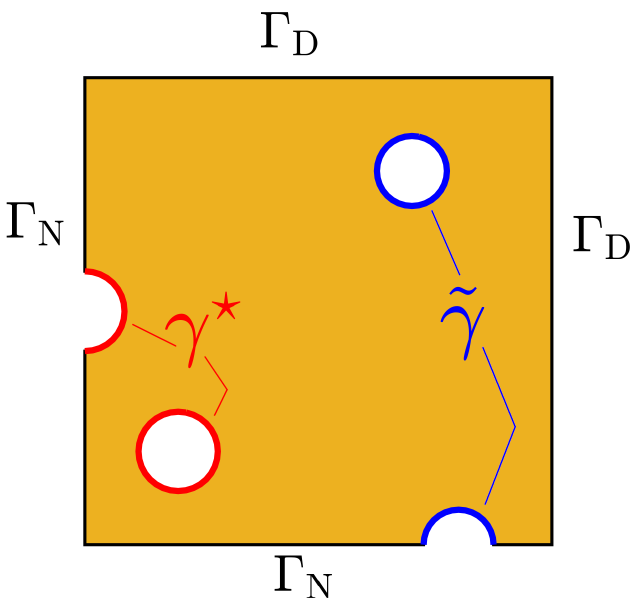}
			\caption{$\Omega$}
			\label{fig:gammargammas}
		\end{subfigure}\hfill
		\begin{subfigure}{.3\textwidth}
			\centering
			\includegraphics[width=0.75\linewidth]{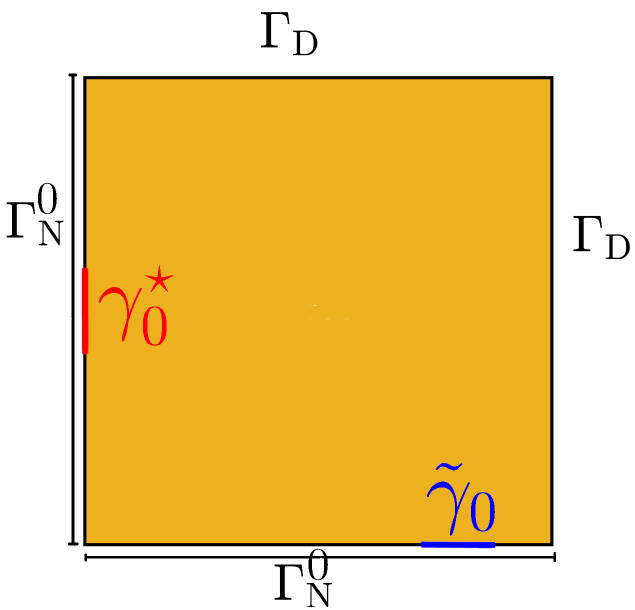}
			\caption{$\Omega_0$}
			\label{fig:nogamma}
		\end{subfigure}\hfill	\begin{subfigure}{.3\textwidth}
		\centering
		\includegraphics[width=0.75\linewidth]{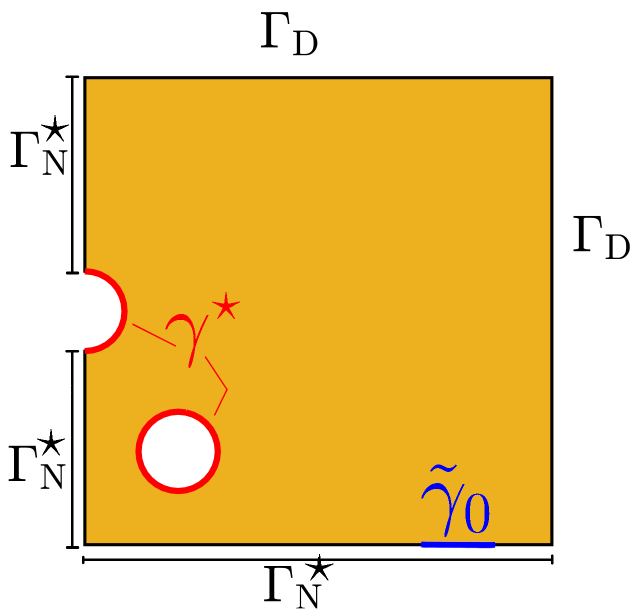}
		\caption{$\Omega^\star$}
		\label{fig:gammas}
	\end{subfigure}
		\caption{Example of a geometry with four negative features, corresponding totally defeatured geometry and partially defeatured geometry retaining only two features. }
		\label{fig:domain}
	\end{figure}

	Let us split the set of the features in two disjoint subsets $\tilde{\mathcal{F}}$ and $\mathcal{F}^\star$, such that $\mathcal{F}=\tilde{\mathcal{F}}\cup \mathcal{F}^\star$ and $\tilde{\mathcal{F}}\cap \mathcal{F}^\star=\emptyset$, and let us define a \emph{partially defeatured geometry as} 
 \begin{equation}
\Omegas:=
\text{int}\big(\overline{\Omega}\cup\lbrace\overline F,~F \in \tilde{\mathcal{F}}\rbrace\big)
\end{equation} 
The partially defeatured geometry can hence be seen either as a simplification of $\Omega$, in which the features in $\tilde{\mathcal{F}}$ are filled with material and the ones in $\mathcal{F}^\star$ are retained, or as a geometrical refinement of $\Omega_0$, in which the features in  $\Fs$ are \emph{included} and the ones in $\Fr$ are \emph{neglected}. To clarify this second nomenclature, that is often used throughout the paper, we refer to Figure \ref{fig:domain}. In $\Omega$ (Figure~\ref{fig:gammargammas}) we say that all the features are \emph{included}, both the ones in $\Fs$ and in $\Fr$. We set 
\begin{equation}
\gamma^\star=\bigcup_{F \in \mathcal{F}^\star}(\partial F\setminus \overline{\gamma^\star_0}), \quad  \tilde{\gamma}:=\bigcup_{F \in \tilde{\mathcal{F}}}(\partial F\setminus \overline{\tilde{\gamma}_0}),\label{gamma_star}
\end{equation} so that $\gamma=\gs\cup\gr$.
In $\Omega_0$ (Figure~\ref{fig:nogamma}), instead, all the features are \emph{neglected}. We introduce
\begin{equation}
{\gamma}^\star_0:=\bigcup_{F \in \mathcal{F}^\star}(\partial F \setminus \overline{\Gamma_{\mathrm{N}}}), \quad \tilde{\gamma}_0:=\bigcup_{F \in  \tilde{\mathcal{F}}}(\partial F \setminus \overline{\Gamma_{\mathrm{N}}})\label{gamma_tilde},
\end{equation} so that $\gamma_0=\gamma_0^\star\cup\gr_0$.  The partially defeatured geometry $\Omegas$ is reported in Figure~\ref{fig:gammas}: only the features in $\Fs$ are \emph{included} whereas the features in $\Fr$ are \emph{neglected}. We here set
$\Gamma_{\mathrm{N}}^\star=\Gamma_{\mathrm{N}}^0\setminus \gamma_0^\star$.

We are now interested in defining the \emph{partially defeatured problem}, which is the problem set on our partial defeatured geometry $\Omegas$. In particular, we consider
	\begin{equation}
	\begin{cases}
	-\Delta \us =f & \text{in}~ \Omega^\star\\
	\us=g_\mathrm{D} &\text{on }  \Gamma_{\mathrm{D}}\\
	\nabla \us \cdot \bm{n}^\star=g_\mathrm{N}& \text{on }\Gamma_{\mathrm{N}}^\star\\
	\nabla \us \cdot \bm{n}^\star=g& \text{on }\gamma^\star
	\end{cases} \label{eq:strong_omega_star}
	\end{equation}
	where $\mathbf{{n}^\star}$ is the outward unitary normal of $\Omega^\star$. 	The variational formulation of Problem \eqref{eq:strong_omega_star} reads: \emph{find $\us \in H^1_{g_\mathrm{D},\Gamma_{\mathrm{D}}}(\Omega^\star)$ that satisfies, $\forall v \in H^1_{0,\Gamma_{\mathrm{D}}}(\Omega^\star)$}
		\begin{equation}
	(\nabla \us,\nabla v)_{\Omegas}=(f,v)_{\Omegas}+\langle g,v\rangle_{\gs}+\langle g_\mathrm{N},v\rangle_{\Gamma_{\mathrm{N}}^\star}. \label{eq:weak_omega_star}
	\end{equation}	\vspace{-0.5cm}
	
Let us now introduce the so called \emph{active mesh}, defined as the set of elements in $\mathcal{T}_h^0$ which have a non-empty intersection with the partially defeatured geometry, i.e.
$$\mathcal{A}_h:=\lbrace K \in \mathcal{T}_h^0:~K \cap \Omega^\star \neq \emptyset\rbrace,$$  and let $\OmegaT=\bigcup_{K \in \mathcal{A}_h}\overline{K}$ be the domain covered by the elements in $\mathcal{A}_h$. Let $h_K=\mathrm{diam}(K)$ and $h=\max_{K \in \mathcal{A}_h}h_K$. For each $K \in \mathcal{A}_h$ we define $\gamma_{K}^\star= K\cap\gamma^\star$, and we introduce the set of \emph{cut elements}
$$\mathcal{G}_h=\lbrace K \in \A{h}: |\gamma_{K}^\star| \neq 0\rbrace.
$$

Let us now introduce the set \begin{equation}Q_h(\OmegaT)=\mathcal{P}_\polydeg(\mathcal{A}_h):=\left\lbrace q_h\in L^2(\OmegaT):~\tr{q_h}{K}\in \mathcal{P}_{\polydeg}(K),~\forall K \in \mathcal{A}_h \right\rbrace,\label{eq:Qh}
\end{equation} and
$$V_h^0(\OmegaT)=\lbrace v \in \mathcal{C}^0(\overline{\OmegaT})\cap Q_h(\OmegaT): \tr{v}{\Gamma_{\mathrm{D}}} =0\rbrace,
$$
$$V_h(\OmegaT)=\lbrace v \in \mathcal{C}^0(\overline{\OmegaT})\cap Q_h(\OmegaT): \tr{v}{\Gamma_{\mathrm{D}}} =g_\mathrm{D}\rbrace,$$
and let us consider the discrete problem: \emph{find $u_h^\star \in V_h(\OmegaT)$ such that}
\begin{equation}
(\nabla u_h^\star,\nabla v_h)_{\Omega^\star}=(f,v_h)_{\Omega^\star}+\langle g,v_h\rangle_{\gs}+\langle g_\mathrm{N},v_h\rangle_{\Gamma_{\mathrm{N}}^\star}\quad \forall v_h \in V_h^0(\OmegaT), \label{eq:prob_num}
\end{equation}
where, for the sake of simplicity, we assume $g$ to be a piecewise linear polynomial on the partition of $\gs$ induced by $\mathcal{A}_h$.

\begin{rem}
	{Let us remark that Equation \eqref{eq:prob_num} is defined on the union of the active elements $\OmegaT$, but the integrals are computed on the partially defeatured domain $\Omegas \subseteq \OmegaT$. If some elements have a very small intersection with $\Omegas$, this may lead to a poorly conditioning in the linear system. It is worth underlining, however, that the Neumann boundary condition prescribed on $\gs$ is imposed weakly (it is a natural boundary condition for the primal solution), so the stability is not affected. Conditioning issues can be mitigated by \textit{ghost penalty} techniques \cite{BURMAN_GHOST}, at the price of modifying the bilinear form in \eqref{eq:prob_num}. Hence, for the sake of simplicity in the developments that follow, we do not consider ghost penalty terms in the primal problem. }
\end{rem}

\section{Error estimator via equilibrated flux reconstruction}\label{sec:flux}

Our aim is to control the energy norm of the error $||\nabla(u-\uhs)||_\Omega$, where $u$ is the solution of \eqref{eq:weak_omega} and $\uhs$ is the solution of \eqref{eq:prob_num}. However, we want to achieve this without solving Problem \eqref{eq:weak_omega}. For this reason we propose an \emph{a posteriori} error estimator defined uniquely on $\Omega^\star$ and based on an equilibrated flux reconstructed from $\uhs$. This estimator allows both for standard mesh refinement and for geometric refinement, which consists in choosing, at each step, which features need to be included into the geometry in order to significantly increase the accuracy of the solution. It is directly adapted from \cite{BCGVV} in order to handle the case of trimmed elements. Indeed, starting from the totally defeatured geometry $\Omega_0$, we aim at gradually including new features without remeshing the domain, but only performing local refinements of the initial non-conforming mesh $\mathcal{T}_h^0$.

 In this section we recall the concept of equilibrated flux reconstruction and we devise a strategy to approximate it on a mesh trimmed by the features in $\mathcal{F}^\star$. The \emph{a posteriori} error analysis based on this flux reconstruction will follow in Section \ref{sec:aposteriori}.\vspace{0.2em}

	Let $\bm{\sigma}\in \Hdiv{\Omega^\star}$ be such that 
		\begin{equation}
	\begin{cases}
	\nabla \cdot \bm{\sigma}=f & \text{ in } \Omega^\star\\
	\bm{\sigma} \cdot\bm{n}=-g_\mathrm{N} & \text{ on }  \Gamma_{\mathrm{N}}^\star\\
	\bm{\sigma} \cdot\bm{n}=-g & \text{ on } \gs.
	\end{cases}\label{eq:flux_prop}
	\end{equation}
	Then, if $\us \in H^1_{g_\mathrm{D},\Gamma_{\mathrm{D}}}(\Omega^\star)$ is the solution of \eqref{eq:weak_omega_star}, the theorem of Prager and Synge (see \cite{Braess,PragerSynge}) states that, for $v \in H^1_{g_\mathrm{D},\Gamma_{\mathrm{D}}}(\Omegas)$,
	\begin{equation}
	||\nabla(\us-v)||_\Omegas^2+||\nabla \us +\bm{\sigma}||_\Omegas^2=||\nabla v +\bm{\sigma}||_\Omegas^2.\label{PragerSynge}
	\end{equation}
	Choosing $v=\uhs$ in $\eqref{PragerSynge}$ we obtain 
	$$||\nabla(\us-\uhs)||_\Omegas^2\leq||\nabla \uhs +\bm{\sigma}||_\Omegas^2,$$
	which means that the difference between the numerical flux and a flux $\bm{\sigma}$ satisfying \eqref{eq:flux_prop} provides a sharp upper bound for the numerical error, with unitary reliability constant.
	
	Equilibrated fluxes are about the cheap construction of such a $\bm{\sigma}$ at a discrete level.
In order to build $\flux$, we resort to a local equilibration procedure directly adapted from \cite{BraessScho2008,EV2015}. Under this approach, local equilibrated fluxes are built on patches of elements sharing a node, avoiding to solve a global optimization problem and making the method well-suited for parallel implementation. 

However, some extra care is needed in order to account for the non-conformity of the mesh with respect to $\gs$. Indeed, the Neumann boundary condition is essential for the flux, and it can be imposed strongly on $\Gamma_{\mathrm{N}}^\star$ in the definition of the discrete space, since we have assumed that the mesh is conforming to this portion of the boundary. The Neumann boundary condition on $\gs$ has instead to be imposed weakly: we hence decide to work in {an unfitted finite element} framework, adopting a Nitsche method as proposed in \cite{puppi2021cut} for the case of the mixed formulation of the Darcy problem. For a flux recovery in the case of {an unfitted Dirichlet boundary} we refer the reader to \cite{HeCapatina}, which resorts however to a different local equilibration procedure.

\begin{rem} Let us remark that the non conformity of the mesh with respect to $\gr$ is instead not relevant when building $\flux$, since the equilibrated flux is reconstructed on the partially defeatured geometry $\Omegas$, which is blind to the features in $\tilde{\mathcal{F}}$.
	
\end{rem}

In the following, we denote by $\N{h}$ and $\mathscr{E}_h$ respectively the set of vertices and faces (edges if $d=2$) in the active mesh $\A{h}$. In particular the set of vertices is decomposed into $\mathscr{N}_h^{\mathrm{ext}}$ and $\mathscr{N}_h^{\mathrm{int}}$, corresponding respectively to vertices lying on $\Gamma_{\mathrm{N}}^\star$ and inside $\Omega_0$. Let us remark that the nodes belonging to $\mathcal{A}_h$ but lying inside a feature will still belong to $\mathscr{N}_h^\mathrm{int}$.

\begin{figure}
	\centering
		\includegraphics[width=0.4\linewidth]{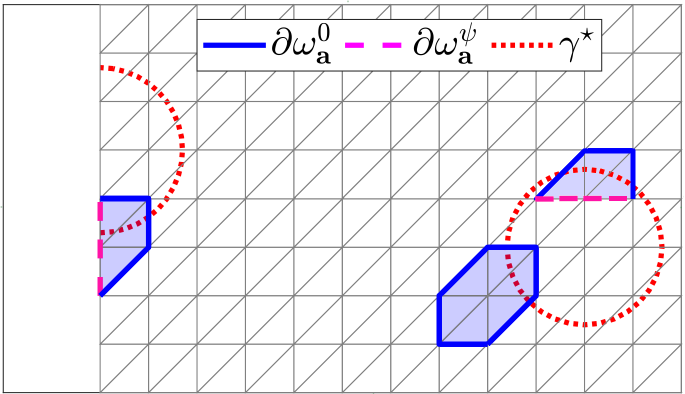}\hspace{0.5cm}%
			\includegraphics[width=0.4\linewidth]{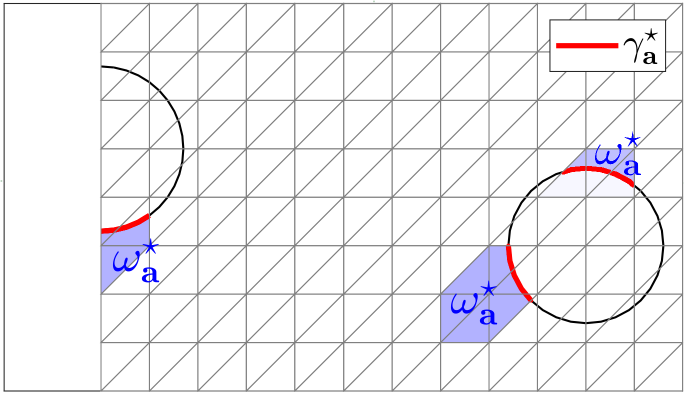}
			
	\caption{Nomenclature for patch boundaries}
	\label{fig:patch_boundary}
\end{figure}

Let us consider a vertex $\bm{a}\in \mathscr{N}_h$ and let us denote by $\wa$ the open patch of elements of $\mathcal{A}_h$ sharing that node, and by $\psi_{\bm{a}}$ the \emph{hat} function in $Q_h(\OmegaT)\cap H^1(\OmegaT)$ taking value 1 in vertex $\bm{a}$ and 0 in all the other vertices. Let $\partwa{}$ be the boundary of $\wa$, and let us define (see Figure~\ref{fig:patch_boundary})
$\partwa{0}\subseteq \partwa{}$ as
$$\partwa{0}=\lbrace \bm{x} \in \partial \wa~:\psi_{\bm{a}}(\bm{x})=0\rbrace,$$ and $\partial \omega_{\bm{a}}^\psi=\partial \wa \setminus {\partial \omega_{\bm{a}}^0}$.
Let us now introduce the space 
$$\bm{M}_h(\wa):=\left\lbrace  \bm{v}_h\in \Hdiv{\wa}:~\tr{\bm{v}_h}{K} \in [\mathcal{P}_\polydeg(K)]^d+\bm{x}\mathcal{P}_\polydeg(K),~\forall K \in \wa\right\rbrace,$$ 
which is the Raviart-Thomas finite element space of order $\polydeg$ on $\wa$, and let
$$\bm{M}_h^{\bm{a},0}=\left\lbrace \bm{v}_h \in \bm{M}_h(\wa):~\bm{v}_h\cdot \bm{n}_{\wa}=0 \text{ on }\partwa{0}\right\rbrace$$ and $Q_h(\wa)$ be the restriction of $Q_h(\OmegaT)$ to $\wa$. We then define
\begin{align}
\Mha&:=\begin{cases}\bm{M}_h^{\bm{a},0}&\text{ if } \bm{a} \in \mathscr{N}_h^\mathrm{int}\\
\big\lbrace \bm{v}_h \in \bm{M}_h(\wa):~\bm{v}_h\cdot \bm{n}_{\wa}=0 \text{ on } \partwa{0} \\\hspace{2.7cm} \bm{v}_h\cdot \bm{n}_{\wa}=-\psia g_\mathrm{N} \text{ on }\partwa{\psi}\cap \Gamma_{\mathrm{N}}^0\rbrace& \text{ if } \bm{a} \in \N{h}^{\mathrm{ext}} ,
\end{cases}\label{eq:Mha}
\end{align}
\begin{align}
\Qha&:=\begin{cases}\big\lbrace q_h \in Q_h(\wa):(q_h,1)_{\wa}=0 \big\rbrace & \text{ if } \bm{a} \in \N{h}^{\mathrm{int}} \text{ or } \bm{a} \in \text{int}(\overline{\Gamma_{\mathrm{N}}^0})\\[0.4em] Q_h(\wa) &\text{ if } \bm{a} \in \N{h}^{\mathrm{ext}} \text{ and } \bm{a} \notin \text{int}(\overline{\Gamma_{\mathrm{N}}^0}),
\end{cases}\label{eq:Qha}
\end{align}
where 
$\bm{n}_\wa$ is the outward unit normal of $\wa$.

For the sake of compactness, let us set (see Figure~\ref{fig:patch_boundary}) $$h_{\bm{a}}=\max_{K \subset \wa}h_K, \quad \omega_{\bm{a}}^\star=\wa \cap \Omega^\star, \quad \gamma_{\bm{a}}^\star=\gs \cap \wa.$$
Let us then introduce the bilinear forms $m_{\bm{a}}:\Mha\times\Mha \rightarrow \mathbb{R}$ and ${b_{\bm{a}}}:\Mha\times\Qha \rightarrow \mathbb{R}$ such that
\begin{align}
& m_{\bm{a}}(\bm{\rho}_h,\bm{v}_h):=(\bm{\rho}_h,\bm{v}_h)_{\was}+\frac{1}{h_{\bm{a}}}\langle\bm{\rho}_h \cdot \bm{n},\bm{v}_h \cdot \bm{n}\rangle_{\gamma_{\bm{a}}^\star} \label{eq:ma}\\
&{b_{\bm{a}}}(\bm{v}_h,q_h):=(q_h,\nabla \cdot \bm{v}_h)_{\was}-{\langle q_h,\bm{v}_h \cdot \bm{n} \rangle_{\gamma_{\bm{a}}^\star}}\label{eq:ba}
\end{align}
and the linear operators $L_{\bm{a}}:\Mha \rightarrow \mathbb{R}$ and ${R_{\bm{a}}}:\Qha \rightarrow \mathbb{R}$ such that
\begin{align}
& L_{\bm{a}}(\bm{v_h}):=-(\psi_{\bm{a}}\nabla u_h^\star,\bm{v}_h)_{\was }-\frac{1}{h_{\bm{a}}}\langle \psi_{\bm{a}}g,\bm{v}_h \cdot \bm{n}\rangle_{\gamma_{\bm{a}}^\star}\label{Lh_nit}\\
& {R_{\bm{a}}}(q_h):=(\psi_{\bm{a}}f-\nabla \psi_{\bm{a}}\cdot \nabla u_h^\star,q_h)_{\was}+{\langle \psi_{\bm{a}}g,q_h\rangle_{\gamma_{\bm{a}}^\star}}\label{Rh_nit},
\end{align}
On each patch we then solve a problem in the form: \emph{find $(\fluxa,\lam)\in \Mha \times \Qha$ such that}
\begin{equation}
\begin{cases}
m_{\bm{a}}(\fluxa,\bm{v}_h)-{b_{\bm{a}}}(\bm{v}_h,\lam)=L_{\bm{a}}(\bm{v}_h) & \forall \bm{v}_h \in \Mha\\
\hspace{2.2cm}{b_{\bm{a}}}(\bm{\fluxa},q_h)={R_{\bm{a}}}(q_h)  & \forall q_h \in \Qha,\label{eq:prob_eq_flux_nit}
\end{cases}
\end{equation}
Finally, the flux reconstruction is built by summing all the local contributions, i.e. 
$$\flux=\sum_{\bm{a}\in \mathscr{N}_h}\fluxa.$$

If $\Fs=\emptyset$, i.e. none of the features in $\mathcal{F}$ is included into the geometry, then it is possible to prove that, given the regularity assumptions on $f$, the flux $\flux$ reconstructed from $\uhs=u_h^0$ is such that $||\nabla \cdot \flux -f||_K=0$, $\forall K \in\mathcal{A}_h= \mathcal{T}_h^0$ (see, e.g., \cite[Lemma 3.5]{EV2015}). 

When $\Fs\neq \emptyset$, the flux reconstructed from $\uhs$ is instead a {weakly equilibrated flux reconstruction}, due to the error introduced by the weak imposition of the Neumann boundary condition on $\gs$. Indeed $||\flux \cdot \bm{n}+g||_{\gs}$ will not be zero, and also the mass balance may be polluted, ending up in some elements for which $||\nabla \cdot \flux-f||_{K \cap \Omegas}\neq 0$. For any $\bm{a} \in \mathscr{N}_h$ we denote by $\overline{q_h}^\wa$ the average of $q_h \in Q_h(\wa)$ on $\wa$. Setting $c_{\bm{a}}=-\overline{q_h}^\wa$, we have that $q_h+c_{\bm{a}}\in \Qha$. The second equation in \eqref{eq:prob_eq_flux_nit} can hence be rewritten as
\begin{align}
(\nabla \cdot \fluxa,q_h+c_{\bm{a}} )_{\was}-\langle\fluxa \cdot \bm{n}+\psia g,q_h+c_{\bm{a}} \rangle_{\gamma_{\bm{a}}^\star}=(\psia f-\nabla \psi_{\bm{a}}\cdot \nabla u_h^\star,q_h+c_{\bm{a}})_{\was} .\label{eq:div_p1}
\end{align}
According to \eqref{eq:prob_num}, we have that
\begin{equation}
(\nabla \psia\cdot\nabla u_h^\star,1)_{\was}=(\psia f,1)_{\was}+\langle \psia g,1\rangle_{\gamma_{\bm{a}}^\star}
+\langle \psia g_\mathrm{N},1\rangle_{\partial\omega_{\bm{a}}^\psi\cap \Gamma_{\mathrm{N}}^\star},\label{eq:hat_orth}
\end{equation}
while, by the divergence theorem, 
\begin{equation}(\nabla \cdot \fluxa,1)_{\was}=\langle\fluxa \cdot \bm{n},1\rangle_{\gamma_{\bm{a}}^\star}+\langle\fluxa \cdot \bm{n}_\wa,1\rangle_{\partial\omega_{\bm{a}}^\psi\cap \Gamma_{\mathrm{N}}^\star}\label{eq:div_theo}.
\end{equation}
Exploiting the fact that the Neumann boundary condition is imposed in strong form on the whole $ \partwa{\psi}\cap \Gamma_{\mathrm{N}}^0$
we can rewrite \eqref{eq:div_p1} as
\begin{align}
(\nabla \cdot \fluxa,q_h )_{\was}-\langle\fluxa \cdot \bm{n}+\psia g,q_h \rangle_{\gamma_{\bm{a}}^\star}=(\psia f-\nabla \psi_{\bm{a}}\cdot \nabla u_h^\star,q_h)_{\was}\quad \forall q_h \in Q_h(\wa). \label{eq:div_intermediate}
\end{align}
Let us consider $K \in \mathcal{A}_h$ and let us denote its nodes by $\mathscr{N}_h(K)$.
As the polynomials in $Q_h(\wa)$ are discontinuous, then \eqref{eq:div_intermediate} holds also for $q_h \in \mathcal{P}_\polydeg(K)$.
We hence have that $\forall  q_h \in \mathcal{P}_\polydeg(K)$
\begin{equation}
(\nabla \cdot \flux,q_h)_{K\cap \Omegas}
=\sum_{\bm{a} \in \mathscr{N}_h(K)}(\nabla \cdot \fluxa,q_h)_{K\cap \Omegas}=(f,q_h)_{K \cap \Omegas}+\langle\flux \cdot \bm{n}+g ,q_h\rangle_{\gamma_{K}^\star} ,\label{eq:div_p1_sum}
\end{equation}
where we have exploited 
the fact that $\psia$ is a partition of unity, and hence $\sum_{\bm{a} \in \mathscr{N}_h(K)}\nabla \psia=0$. This means that, since we have assumed that $f \in Q_h(\Omega_0)$,
\begin{equation}
||\nabla \cdot \flux-f||_K=0\quad \forall K \in \mathcal{A}_h \setminus \mathcal{G}_h. \label{eq:zero_div}
\end{equation}
{This means that equilibration holds true on every uncut element, failing instead on cut elements.}

\begin{rem} \label{rem:p0}
	The weakly equilibrated flux reconstruction could be built also by considering asymmetric local patch problems instead of \eqref{eq:prob_eq_flux_nit}, essentially by omitting the integrals over $\gamma_{\bm{a}}^\star$ in the mass balance equation \cite{puppi2021cut}. However, as shown in Appendix \ref{appendix1}, this would imply a mass imbalance on a larger set of elements, namely not only in the cut elements, but in all the the elements belonging to a cut patch ($\wa$ such that $|\gamma_{\bm{a}}^\star|\neq 0$). For this reason the symmetric formulation \eqref{eq:prob_eq_flux_nit} is preferred.
\end{rem}

\begin{rem}\label{stab}
	{The unfitted finite element strategy proposed in \cite{puppi2021cut} and adopted herein is based on CutFEM, which is known to require stabilization}
	in presence of very small cuts {when essential boundary conditions are imposed on the unfitted boundary} \cite{BURMAN_GHOST} . Let $v :\OmegaT\rightarrow \mathbb{R}$ be smooth enough and let us define the \emph{jump} as
	$$\llbracket v\rrbracket_e:=\tr{v}{K_1}-\tr{v}{K_2},$$ with $e=\partial K_1\cap \partial K_2$. Let us then introduce for $\bm{a} \in \mathscr{N}_h$,
	$$\mathscr{E}_{h,\bm{a}}:=\lbrace e\in \mathscr{E}_h \text{ such that  }~ \exists K\in \mathcal{A}_h,~K\subset \wa:~e\subset \partial K \rbrace,$$ and let
 $\mathscr{E}_{h,\bm{a}}^{\mathrm{int}}:=\lbrace e\in \mathscr{E}_{h,\bm{a}} \text{ such that  }~ e\nsubseteq\partial \wa\rbrace.$ Finally we set
$$\mathscr{E}_{h,\bm{a}}^{\mathcal{G}}:=\lbrace e\in \mathscr{E}_{h,\bm{a}}^{\mathrm{int}}\text{ such that  }~ \exists K\in \mathcal{G}_h:~e\subset \partial K\rbrace$$ 
and we define the bilinear forms (\cite{BURMAN_GHOST,BurmanGeometry}) $\bm{j}_h:\Mha \times \Mha\rightarrow \mathbb{R}$, $j_h:\Qha\times \Qha\rightarrow \mathbb{R}$ such that 
	\begin{align}
	&\bm{j}_{\bm{a}}(\bm{\rho}_h,\bm{v}_h):=\beta_1\sum_{e \in \mathscr{E}_{h,\bm{a}}^{\mathcal{G}}}\Big[h_{\bm{a}}(\llbracket\bm{\rho}_h\rrbracket_e, \llbracket\bm{v}_h\rrbracket_e)_e+\sum_{j=1}^rh_{\bm{a}}^{2j+1}(\llbracket D_{\bm{n}_e}^j\bm{\rho}_h \rrbracket_e,\llbracket D_{\bm{n}_e}^j\bm{v}_h \rrbracket_e)_e\Big]\label{stab1}\\
	&j_{\bm{a}}(\mu_h,q_h):=\beta_2\sum_{e \in \mathscr{E}_{h,\bm{a}}^{\mathcal{G}}}\Big[\frac{1}{h_{\bm{a}}}(\llbracket\mu_h\rrbracket_e, \llbracket q_h\rrbracket_e)_e+\sum_{j=1}^r h_{\bm{a}}^{2j-1}(\llbracket D_{\bm{n}_e}^j \mu_h \rrbracket_e,\llbracket D_{\bm{n}_e}^j q_h \rrbracket_e)_e\Big],\label{stab2}
	\end{align}
	where $\beta_1,\beta_2>0$ and $D_{\bm{n}_e}^j$ is the normal derivative of order $j$, with $\bm{n}_e$ denoting the unit normal vector to $e$, whose orientation is fixed.
	The stabilized version of Problem \eqref{eq:prob_eq_flux_nit} then reads as: \emph{find $(\fluxa,\lam)\in \Mha \times \Qha$ such that}
	\begin{equation}
	\begin{cases}
	m_{\bm{a}}(\fluxa,\bm{v}_h)-{b_{\bm{a}}}(\bm{v}_h,\lam)+\bm{j}_{\bm{a}}(\flux,\bm{v}_h)=L_{\bm{a}}(\bm{v}_h) & \forall \bm{v}_h \in \Mha\\
	\hspace{2.2cm}{b_{\bm{a}}}(\bm{\fluxa},q_h)-j_{\bm{a}}(\lam,q_h)={R_{\bm{a}}}(q_h)  & \forall q_h \in \Qha.\label{eq:prob_eq_flux_nit_stab}
	\end{cases}
	\end{equation}
The scalar stabilization term introduced in the second equation is however likely to further pollute mass balance, not only on cut elements but on all the cut patches. Using a projection-based stabilization only on badly cut patches (see for example \cite{puppi2021cut}) would allow to reduce the number of elements on which this additional error is introduced, but it might also affect the monotonicity of the estimator proposed in the next section. Indeed, assuming to use the estimator in an adaptive procedure, badly cut elements are very likely to appear and disappear as the mesh changes, introducing extra contributions to the estimator only at certain steps.
	 Divergence preserving techniques have been proposed in \cite{Zahedi_divFREE}, requiring however the evaluation of higher order derivatives, and hence increasing the computational effort required to compute $\flux$. Let us recall that $\flux$ is not the primal variable of our problem: it is computed only with the aim of building the estimator, so it is of interest to keep its computation as cheap as possible. For this reason, badly cut elements could also be simply neglected, allowing to avoid the case in which an element that has barely an active part, strongly affects the total estimator only due to stabilization. Further remarks on this aspect can be found in Section \ref{sec:num_exp}.
	\end{rem}
	\begin{rem}\label{rem:pos}
	In the present work we decided to focus on the negative Neumann feature case and on the challenges posed by the presence of cut elements when essential boundary conditions need to be weakly imposed in the equilibrated flux reconstruction. Referring to the terminology introduced in \cite{BCV2022}, also positive features could be considered, in the sense of small additions of material at the boundary of the domain. In that case the equilibrated flux reconstruction in the positive feature and in the rest of the domain could be computed separately, provided that proper coupling conditions are defined at the interface. What is interesting to remark is that, as reported in \cite{BCGVV} (see also \cite{BCV2022}), positive features with complex shapes are often treated by introducing a bounding box having the simplest possible shape, meshing the bounding box instead of the actual feature and computing the flux reconstruction on that mesh, which is blind to the feature boundary. Since a positive feature is actually a negative feature from its bounding box perspective, all the analysis that was carried out until now can easily be exploited. In contrast, the a posteriori error analysis for the Dirichlet feature case is more involved, as it requires to work on Sobolev spaces with negative exponent. We refer to \cite{weder2025} for the derivation of the defeaturing component of the estimator in that case.
\end{rem}

\section{\emph{A posteriori} error analysis}\label{sec:aposteriori}
In this section we propose a reliable estimator for the error $||\nabla( u-\uhs)||_\Omega$, based on the (weakly) equilibrated flux defined above.

For $F \in \Fr$ let ${\gr_{F}}:=\gr\cap\partial F$ and $\gr_{{0,F}}=\gr_0\cap F$, so that $\gr=\bigcup_{F \in \Fr}{\gr_{F}}$ and $\gr_0=\bigcup_{F \in \Fr}{\gr_{0,F}}$. Let $\bm{\sigma}\in \Hdiv{\Omegas}$ be a generic flux satisfying \eqref{eq:flux_prop}. Then, for every $F \in \Fr$ we denote by $\tr{\bm{\sigma}}{F}\in \Hdiv{F}$ its restriction to $F$, satisfying 
	\begin{equation}
\begin{cases}
\nabla \cdot\tr{\bm{\sigma}}{F}=f & \text{ in } F\\
\tr{\bm{\sigma}}{F} \cdot\bm{n_0}=-g_0 & \text{ on }  \gr_{0,F}\\
\tr{\bm{\sigma}}{F} \cdot\bm{n_F}=-\bm{\sigma} \cdot\bm{n} & \text{ on } \gr_F.
\end{cases}\label{eq:flux_prop_restricted}
\end{equation} Vectors $\bm{n_F}=-\bm{n}$ and $\bm{n_0}$ are the unitary outward normals of $F$ on $\gr_F$ and $\gr_{0,F}$, respectively. Omitting, by an abuse of notation, the explicit restriction of $\bm{\sigma}$ to $F$, we define for every $F \in \Fr$, the quantity
$$d_F:=g+\bm{\sigma}\cdot \bm{n} \quad \text{on } \gr_{F}, $$
	 which is the error between the Neumann datum $g$ on ${\gr_{F}}$ and the normal trace of $\bm{\sigma}$. Similarly we define
	 $$
	 d_F^h:=g+\flux\cdot \bm{n} \quad \text{on } \gr_{F},$$ where $\flux$ is a (weakly) equilibrated flux reconstructed from $\uhs$ as detailed in Section \ref{sec:flux}. Denoting by $\overline{d_F}^{\gr_F}$ and $\overline{d_F^h}^{\gr_F}$ the average, respectively, of $d_{F}$ and $d_F^h$ over ${\gr_{F}}$ we define
	 \begin{equation}
	 \estim{F}^2:=|{\gr_{F}}|^{\frac{1}{d-1}}\big|\big|d_F^h-\overline{d_F^h}^{\gr_F}\big|\big|_{{\gr_{F}}}^2+c_{{\gr_{F}}}^2|{\gr_{F}}|^{\frac{d}{d-1}}\big|\overline{d_F}^{\gr_F}\big|^2\label{eq:estim_gamma},
	 \end{equation}
	 with
	 $c_{{\gr_{F}}}$ defined according to \eqref{c_omega}. Let us remark that $\overline{d_F}^{\gr_F}$ actually depends only on the choice of the extension of $f$ inside $F$ and on the choice of $g_0$. Indeed, according to \eqref{eq:flux_prop_restricted}, $$\overline{d_F}^{\gr_F}=\frac{1}{|\gr_{F}|}(g+\bm{\sigma}\cdot \bm{n},1)_{\gr_{F}}=\frac{1}{|\gr_{F}|}\big[(g,1)_{\gr_{F}}-(f,1)_{F}-(g_0,1)_{\gr_{0,F}}\big],$$ so that a proper choice of the data allows us to get rid of the second term in $\estim{F}^2$. 
	 
	 For each $K \in \mathcal{A}_h$ we also define
	 	\begin{equation}
	 \estim{\mathrm{div}}^K:=h_K||f-\nabla \cdot\flux||_{K\cap \Omegas}\qquad \estim{\sigma}^K:=||\flux+\nabla \uhs||_{K\cap \Omegas}. 
	 \label{estim_div_and_g}\vspace{-0.2cm}
	 \end{equation}
	 \begin{equation}
	 \estim{\mathrm{g}}^K:=h_K^\frac{1}{2}(||g+\flux\cdot \bm{n}||_{\gamma_{K}^\star}^2\big)^\frac{1}{2}\label{estim_sigma}
	 \end{equation}

	 We now state and prove the proposition establishing our \emph{a posteriori} bound under {two technical assumptions} and a lemma proved in \cite{GuzmanOlshanskii2018}:
	 \begin{assum}\label{assum:hf}
	 	For every $F \in \Fr$ let ${\delta_F}:=\mathrm{diam}(F)$ and ${h_{F}^\mathrm{min}}:=\min\lbrace{h_K:K \in \mathcal{A}_h,~K\cap F\neq \emptyset}\rbrace$. Then we assume that ${\delta_F} \lesssim {h_{F}^\mathrm{min}}$. {In view of the adaptive strategy detailed in Section \ref{sec:adaptive} (see also Remark \ref{rem:assum}), this corresponds to say} that the number of elements covering a non-included feature can not grow indefinitely.
	 \end{assum}
	 {\begin{assum}\label{assum4}
	 	For every $F \in \mathcal{F}^\star$ such that $\gamma^\star \subset \partial F$ let $h_F^{\rm{max}}:=\max\lbrace h_K: K \in \mathcal{G}_h,~|\gamma_K^\star|\neq 0\rbrace$. Then we assume that $h_F^{\rm{max}}\lesssim \delta_F$. In view of the adaptive strategy detailed in Section \ref{sec:adaptive} (see also Remark \ref{rem:assum}), this corresponds to say that an included feature is resolved by the surrounding mesh.
	 \end{assum}
 	\begin{lemma}\label{lem:loctrace}
 		Let $F \in \mathcal{F}^\star$ be an included feature, let $\delta_F:=\mathrm{diam}(F)$, and let
 		$\hat F := \delta_F^{-1}F$ denote the corresponding feature rescaled to unit diameter. Let
 		$h_0>0$ and let $K \in \mathcal{G}_h$ be such that $h_K < h_0$. Then, for every
 		$v \in H^1(K)$,
 		\begin{equation}\label{eq:loc_trace_lem}
 			\|v\|^2_{\gamma_K^\star} \leq C_{\rm tr}\big(h_K^{-1}\|v\|^2_{K} + h_K\|\nabla v\|^2_{K}\big),
 			\qquad
 			C_{\rm tr} := C_{\hat F}\,C_R^2\big(1+2h_0\delta_F^{-1}\big),
 		\end{equation}
 		where $C_R$ depends only on $d$ and on the shape regularity of $\mathcal{T}_h^0$, and
 		$C_{\hat F}$ only on $d$ and on the shape of $\hat F$. In particular, neither
 		constant depends on $h$, on $\delta_F$, on $\mathrm{diam}(\Omega)$, nor on how $\gamma^\star$
 		intersects $K$.
 \end{lemma}}

 {For the sake of completeness, and to explicitly track the constants, we report the proof of Lemma 1 in Appendix \ref{appendix2}. Let us remark that $C_{\mathrm{tr}}$ depends on the choice of $h_0>0$. Under Assumption \ref{assum4}, we can choose any $h_0>h_F^{\rm{max}}$ and obtain a constant that is bounded independently of $h$ and $\delta_F$.}
	  \begin{prop}\label{prop1}
	 	Let $u$ be the solution to \eqref{eq:weak_omega}, $\uhs$ the solution to \eqref{eq:prob_num}, and let $\flux$ be a weakly equilibrated flux reconstructed from $\nabla \uhs$ by solving the patch local problems in \eqref{eq:prob_eq_flux_nit}. Then, under Assumptions \ref{assum:hf} {and \ref{assum4}}, there exist three constants $\alpha_1,\alpha_2,\alpha_3>0$ such that
	 	\begin{align*}
	 	||\nabla(u-u_h^\star)||_{\Omega}\leq
	 \Big(\sum_{K \in \mathcal{G}_h}\alpha_1(\estim{\mathrm{div}}^K)^2\Big)^{\frac{1}{2}}+\Big(\sum_{K \in \mathcal{G}_h}\alpha_2(\estim{\mathrm{g}}^K)^2\Big)^{\frac{1}{2}}+
	 		\Big(\sum_{K \in \mathcal{A}_h}(\estim{\sigma}^K)^2\Big)^\frac{1}{2}+\Big(\sum_{F \in \Fr}\alpha_3\estim{F}^2\Big)^\frac{1}{2}.
	 	\end{align*}
	 	All the constants $\alpha_i$, $i=1,2,3$ may depend on the spatial dimension $d=2,3$, on the {shape of the features} and on the polynomial degree $r$, with $\alpha_1$ and $\alpha_2$ depending also on the shape regularity of the mesh $\T{h}^0$. They are all independent on the feature size, on the mesh size, and on how the feature boundaries cut the mesh.	 
	 \end{prop}
 \begin{proof}
 		Let $v \in H^{1}_{0,\Gamma_{\mathrm{D}}}(\Omega)$. Adding and subtracting $(\flux,\nabla v)_\Omega$, and exploiting \eqref{eq:weak_omega} we have that
 	\begin{align}
 	(\nabla (u-\uhs), \nabla v)_\Omega&=(\nabla u +\flux,\nabla v)_\Omega-(\nabla \uhs+\flux, \nabla v)_\Omega\nonumber\\
 	&=(f-\nabla \cdot \flux,v)_\Omega+\langle g+\flux\cdot \bm{n},v\rangle_{\gs}+\langle g+\flux\cdot \bm{n},v\rangle_{\gr}\nonumber\\&\quad-(\nabla \uhs+\flux, \nabla v)_\Omega\label{eq:step1_frfs}.
 	\end{align}
 		{Let $\mathbf{E}^0:H_{0,\Gamma_{\mathrm{D}}}^{1}(\Omega)\rightarrow H_{0,\Gamma_{\mathrm{D}}}^{1}(\Omega_0)$ be the extension operator defined in \eqref{stein}.} Let $v_0=\mathbf{E}^0(v) \in H_{0,\Gamma_{\mathrm{D}}}^{{1}}(\Omega_0)$, and let us rewrite \eqref{eq:step1_frfs} as:
 			\begin{align}
 		(\nabla (u-\uhs), \nabla v)_\Omega
 		&=(f-\nabla \cdot \flux,v_0)_\Omegas+\langle g+\flux\cdot \bm{n},v_0\rangle_{\gs}-(\nabla \uhs+\flux, \nabla v_0)_\Omega\nonumber\\
 		&\quad+\sum_{F \in \mathcal{\Fr}}\big[\langle g+\flux\cdot \bm{n},v_0\rangle_{{\gr_{F}}}-(f-\nabla \cdot \flux,v_0)_F\big]= \mathrm{I}+\mathrm{II}+\mathrm{III}\label{eq:step2_frfs},
 		\end{align}with
 		\begin{align*}
 		&\mathrm{I}=(f-\nabla \cdot \flux,v_0)_\Omegas+\langle g+\flux\cdot \bm{n},v_0\rangle_{\gs},\\
 		&\mathrm{II}=-(\nabla \uhs+\flux, \nabla v_0)_\Omega,\\
 		&\mathrm{III}=\sum_{F \in \mathcal{\Fr}}\big[\langle g+\flux\cdot \bm{n},v_0\rangle_{{\gr_{F}}}-(f-\nabla \cdot \flux,v_0)_F\big]\nonumber.
 		\end{align*}
 		
{Let us now consider a projection operator $I_h: H^1(\Omega_0)\rightarrow Q_h(\Omega_0)$ such that $(\nabla I_h v_0,\nabla q_h)_K = (\nabla v_0,\nabla q_h)_K$ $\forall q_h \in Q_h(\Omega_0)$ {and $(v_0-I_hv_0,1)_K=0$ $\forall K \in \mathcal{T}_h^0$}}. In particular, $\forall v_0 \in {H^1}(\Omega_0)$ we have 
 \begin{equation}
 \sum_{K \in \mathcal{T}_h^0}h_K^{-2}||v_0-I_hv_0||_{K}^2\leq{\frac{1}{\pi^2}}||\nabla v_0||_{\Omega_0}^2, \qquad \sum_{K \in \mathcal{T}_h^0}||\nabla I_hv_0||_{K}^2{~\leq~}||\nabla v_0||_{\Omega_0}^2,\label{sc_zhang1}
 \end{equation}
and
 \begin{equation}
 \sum_{K \in \mathcal{T}_h^0}||\nabla(v_0-I_hv_0)||_{K}^2 ~{\leq~}||\nabla v_0||_{\Omega_0}^2.\label{sc_zhang2}
 \end{equation}
 {where the first inequality follows from \cite{Payne1960}, {using the convexity of mesh elements and the fact that $v_0-I_hv_0$ has zero average, together with \eqref{sc_zhang2}}}. Let us remark that both \eqref{sc_zhang1} and \eqref{sc_zhang2} refer to the totally defeatured domain $\Omega_0$, so that there is no dependence on the position of the features with respect to the mesh. 
Since $\tr{I_hv_0}{K}\in \mathcal{P}_\polydeg(K)$, $\forall K \in \mathcal{A}_h$, then according to \eqref{eq:div_p1_sum}, 
\begin{align}
(f-\nabla \cdot \flux,I_hv_0)_{K\cap \Omegas}+\langle\flux \cdot \bm{n}+g ,I_hv_0\rangle_{\gamma_{K}^\star}=0.
\end{align}
Going back to I, we can rewrite its expression as
\begin{align*}
\mathrm{I}&=\sum_{K \in \mathcal{A}_h}(f-\nabla \cdot\flux,v_0-I_hv_0)_{K\cap \Omegas}+\sum_{K \in \mathcal{G}_h}\langle g+\flux \cdot \bm{n},v_0-I_hv_0\rangle_{\gamma_K^\star}=\mathrm{I}_a+\mathrm{I}_b
\end{align*}
with
\begin{align*}
&\mathrm{I}_a=\sum_{K \in \mathcal{A}_h}(f-\nabla \cdot\flux,v_0-I_hv_0)_{K\cap \Omegas}\\
&\mathrm{I}_b=\sum_{K \in \mathcal{G}_h}\langle g+\flux \cdot \bm{n},v_0-I_hv_0\rangle_{\gamma_K^\star}
\end{align*}
	Exploiting the Cauchy--Schwarz inequality, \eqref{eq:zero_div}, \eqref{sc_zhang1} and \eqref{stein} we have that
\begin{align}
\mathrm{I}_a&\leq\sum_{K \in \mathcal{G}_h} ||f-\nabla \cdot \flux||_{K\cap \Omega^\star}||v_0-I_hv_0||_{K\cap \Omega^\star}\nonumber \\
&{\leq} \Big(\sum_{K \in \mathcal{G}_h}  h_K^2||f-\nabla \cdot \flux||^2_{K\cap \Omega^\star}\Big)^{\frac{1}{2}}\Big(\sum_{K \in \mathcal{G}_h}h_K^{-2}||v_0-I_hv_0||^2_{K\cap \Omega^\star}\Big)^\frac{1}{2}\nonumber\\
&\leq \Big(\sum_{K \in \mathcal{G}_h} (\estim{\mathrm{div}}^K)^2\Big)^{\frac{1}{2}}\Big(\sum_{K \in \mathcal{T}_h^0}h_K^{-2}||v_0-I_hv_0||^2_{K}\Big)^\frac{1}{2}\nonumber\\
&{\leq \frac{1}{\pi}} \Big(\sum_{K \in \mathcal{G}_h} (\estim{\mathrm{div}}^K)^2\Big)^{\frac{1}{2}}||\nabla v_0||_{\Omega_0}{\leq \frac{C_\mathrm{stein}}{\pi}} \Big(\sum_{K \in \mathcal{G}_h} (\estim{\mathrm{div}}^K)^2\Big)^{\frac{1}{2}}||\nabla v||_{\Omega},\label{eq:err_Ia}
\end{align}
	For $\mathrm{I}_b$ we have instead that, by the Cauchy--Schwarz inequality
\begin{align}
\mathrm{I}_b&\leq\sum_{K \in \mathcal{G}_h}||g+\flux \cdot \bm{n}||_{\gamma_{K}^\star}||v_0-I_hv_0||_{\gamma_{K}^\star}\nonumber\\
&\leq \Big(\sum_{K \in \mathcal{G}_h}  h_K||g+\flux\cdot \bm{n}||^2_{\gamma_{K}^\star}\Big)^{\frac{1}{2}}\Big(\sum_{K \in \mathcal{G}_h}h_K^{-1}||v_0-I_hv_0||^2_{\gamma_{K}^\star}\Big)^\frac{1}{2}\label{flag}
\end{align}

Exploiting again \eqref{sc_zhang1}-\eqref{sc_zhang2} and Lemma~\ref{lem:loctrace}, it follows that
\begin{align*}
\sum_{K \in \mathcal{G}_h}h_K^{-1}||v_0-I_hv_0||_{\gamma_{K}^\star}^2&{\leq C_\mathrm{tr}}\sum_{K \in \mathcal{G}_h}\big(h_K^{-2}||v_0-I_hv_0||_K^2+||\nabla (v_0-I_hv_0)||_K^2\big)\nonumber \\
&\leq {C_\mathrm{tr}}\sum_{K \in \mathcal{T}_h^0}\big(h_K^{-2}||v_0-I_hv_0||_K^2+||\nabla (v_0-I_hv_0)||_K^2\big)\nonumber \\
&{\leq C_\mathrm{tr}\big(\frac{1 }{\pi^2}+1\big)}||\nabla v_0||_{\Omega_0}^2,
\end{align*} 
Substituting in \eqref{flag} and exploiting \eqref{stein} we then have
\begin{align}
\mathrm{I}_b {\leq C_\mathrm{stein} C_\mathrm{tr}\big(\frac{1 }{\pi^2}+1\big)} \Big(\sum_{K \in \mathcal{G}_h} (\estim{\mathrm{g}}^K)^2\Big)^{\frac{1}{2}} ||\nabla v||_{\Omega}.\label{eq:err_Ib}
\end{align}
For what concerns the term II in \eqref{eq:step2_frfs} we simply apply the Cauchy--Schwarz inequality, obtaining
\begin{align}
\mathrm{II}&\leq ||\flux +\nabla \uhs||_\Omega||\nabla v_0||_\Omega\leq ||\flux +\nabla \uhs||_{\Omegas}||\nabla v||_{\Omega} {=}\Big(\sum_{K \in \mathcal{A}_h}||\flux +\nabla \uhs||^2_{K \cap \Omegas}\Big)^\frac{1}{2}||\nabla v||_\Omega,\nonumber \\[-0.8em]
&=\Big(\sum_{K \in \mathcal{A}_h}(\estim{\sigma}^K)^2\Big)^\frac{1}{2}||\nabla v||_\Omega\label{eq:err_II}
\end{align}
where we have also exploited the fact that $\Omega\subset\Omegas$.

Finally, let us consider the term III of \eqref{eq:step2_frfs}. First of all let us observe that, for any $F \in \Fr$ and for any constant $c\in \mathbb{R}$ we have
\begin{equation*}
(f-\nabla \cdot \flux ,c)_F=(\nabla \cdot \bm{\sigma}-\nabla \cdot \flux,c)_F=-\langle\bm{\sigma}\cdot \bm{n}-\flux\cdot \bm{n},c\rangle_{\gr_F},
\end{equation*}where we have used \eqref{eq:flux_prop_restricted}, and the fact that $\bm{n_F}=-\bm{n}$ .
Denoting by $\overline{v_0}^{\gr_F}$ the average of $v_0$ over $\gr_F$ and choosing $c=-\overline{v_0}^{\gr_F}$, we can rewrite III as
\begin{align*}
\mathrm{III}&=\sum_{F \in \mathcal{\Fr}}\big[\langle g+\flux\cdot \bm{n},v_0\rangle_{{\gr_{F}}}-(f-\nabla \cdot \flux,v_0-\overline{v_0}^{\gr_F})_F+\langle(\bm{\sigma}-\flux)\cdot \bm{n},\overline{v_0}^{\gr_F}\rangle_{\gr_F}\big]=\mathrm{III}_a+\mathrm{III}_b,
\end{align*}
with 
\begin{align*}
&\mathrm{III}_a=\sum_{F \in \Fr}(\nabla\cdot \flux-f,v_0-\overline{v_0}^{\gr_F})_F\\
&\mathrm{III}_b=\sum_{F \in \Fr}\big[\langle g+\flux\cdot \bm{n},v_0\rangle_{{\gr_{F}}}+\langle(\bm{\sigma}-\flux)\cdot \bm{n},\overline{v_0}^{\gr_F}\rangle_{\gr_F}\big]
\end{align*} In order to bound $\mathrm{III}_a$ we use the Friedrichs inequality reported in \cite[Lemma A.2]{BCV2022_arxiv}: since $\gr_F\subset \partial F$, $|\gr_F|>0$ and $|\gr_F|^\frac{1}{d-1}\simeq {\delta_F}$, then $\forall v \in H^1(F)$ 
\begin{equation}
||v-\overline{v}^{\gr_F}||_F{\leq C_\mathrm{f}}~ {\delta_F}||\nabla v||_F, \label{Friedrichs}
\end{equation} where $C_\mathrm{f}$ is independent of the size of $F$, {but can depend on its shape}. We hence obtain, under Assumption \ref{assum:hf} and using the Cauchy--Schwartz inequality, \eqref{Friedrichs}, \eqref{eq:zero_div} and \eqref{stein}, that
\begin{align}
\mathrm{III}_a&\leq \sum_{F \in \mathcal{\Fr}}||\nabla \cdot \flux-f||_F||v_0-\overline{v_0}^{\gr_F}||_F\lesssim \sum_{F \in \mathcal{\Fr}}||\nabla \cdot \flux-f||_F\ {\delta_F}||\nabla v_0||_F\nonumber \\
&\leq\Big(\sum_{F \in \Fr}{\delta_F^2}||\nabla \cdot \flux-f||_F^2\Big)^\frac{1}{2}\Big(\sum_{F \in \Fr}||\nabla v_0||_F^2\Big)^\frac{1}{2}\nonumber\\
&\lesssim\Big(\sum_{F \in \Fr}(h_F^{\min})^2\sum_{K \in \mathcal{G}_h}||\nabla \cdot \flux-f||_{F\cap K}^2\Big)^\frac{1}{2}\Big(\sum_{F \in \Fr}||\nabla v_0||_F^2\Big)^\frac{1}{2}\nonumber\\
&\leq\Big(\sum_{F \in \Fr}\sum_{K \in \mathcal{G}_h}h_K^2||\nabla \cdot \flux-f||_{F\cap K}^2\Big)^\frac{1}{2}\Big(\sum_{F \in \Fr}||\nabla v_0||_F^2\Big)^\frac{1}{2}\nonumber\\
&\leq\Big(\sum_{K \in \mathcal{G}_h}h_K^2||\nabla \cdot \flux-f||_{K\cap\Omegas}^2\Big)^\frac{1}{2}||\nabla v_0||_{\Omega_0}\lesssim\Big(\sum_{K \in \mathcal{G}_h}(\estim{\mathrm{div}}^K)^2\Big)^\frac{1}{2}||\nabla v||_{\Omega},\label{eq:err_IIIa}
\end{align}
with the last hidden constant depending on $C_\mathrm{f}$, $C_\mathrm{stein}$ and the constant in Assumption \ref{assum:hf}.
For what concerns $\mathrm{III}_b$ let us remark that $d_F=d_F^h+(\bm{\sigma}-\flux)\cdot \bm{n}$ on $\gr_F$, $\forall F \in \Fr$. We hence have that
\begin{align*}
\mathrm{III}_b&=\sum_{F \in \Fr}\big[\langle d_F^h,v_0\rangle_{\gr_F}+\langle(\bm{\sigma}-\flux)\cdot \bm{n},\overline{v_0}^{\gr_F}\rangle_{{\gr_{F}}}\big]=\sum_{F \in \Fr}\big[\langle d_F^h,v_0\rangle_{\gr_F}+\langle d_F-d_F^h,\overline{v_0}^{\gr_F}\rangle_{{\gr_{F}}}\big]\\
&=\sum_{F \in \Fr}\big[\langle d_F^h,v_0-\overline{v_0}^{\gr_F}\rangle_{\gr_F}+\langle d_F,\overline{v_0}^{\gr_F}\rangle_{{\gr_{F}}}\big]=\sum_{F \in \Fr}\big[\langle d_F^h-\overline{d_F^h}^{\gr_F},v_0-\overline{v_0}^{\gr_F}\rangle_{\gr_F}+\langle \overline{d_F}^{\gr_F},v_0\rangle_{{\gr_{F}}}\big].
\end{align*} Once $\mathrm{III}_b$ is rewritten in this form, referring the reader to \cite[Theorem 4.3]{BCV2022} it is possible to prove that
 	 	\begin{equation}
 	\mathrm{III}_b\lesssim\Big(\sum_{F \in \Fr} \estim{F}^2\Big)^\frac{1}{2}||\nabla v ||_\Omega\label{eq:err_IIIb},
 	\end{equation}
 	with $\estim{F}^2$ defined as in \eqref{eq:estim_gamma} and the hidden constant depending on the Poincaré constant on $\tilde{\gamma}_F$ (see \cite[Appendix A.1]{BCV2022}) and a trace constant on $\Omega$ .
 	
 	Going finally back to \eqref{eq:step2_frfs} with $v=e:=u-\uhs$ and using \eqref{eq:err_Ia},  \eqref{eq:err_Ib}, \eqref{eq:err_II},  \eqref{eq:err_IIIa},  \eqref{eq:err_IIIb} we can conclude that there exist $\alpha_1, \alpha_2, \alpha_3 >0$ such that
 		\begin{align*}
||\nabla e||_{\Omega}^2&\leq \Big(\sum_{K \in \mathcal{G}_h} \alpha_1 (\estim{\mathrm{div}}^K)^2\Big)^\frac{1}{2}||\nabla e||_{\Omega}+\Big(\sum_{K \in \mathcal{G}_h}\alpha_2 (\estim{\mathrm{g}}^K)^2\Big)^\frac{1}{2}||\nabla e||_{\Omega}+\Big(\sum_{K\in \mathcal{A}_h} (\estim{\sigma}^K)^2\Big)^\frac{1}{2}||\nabla e||_{\Omega}\\&+\Big(\sum_{F \in \Fr} \alpha_3\estim{F}^2\Big)^\frac{1}{2}||\nabla e||_{\Omega}.
 		\end{align*}
 		Simplifying on both sides yields the thesis.
 	\end{proof}
In the following, we will refer to 
\begin{align*}\estim{\mathrm{num}}&:=\Big(\sum_{K \in \mathcal{G}_h}  \alpha_1(\estim{\mathrm{div}}^K)^2\Big)^\frac{1}{2}+\Big(\sum_{K \in \mathcal{G}_h}\alpha_2 (\estim{\mathrm{g}}^K)^2\Big)^\frac{1}{2}+\Big(\sum_{K\in \mathcal{A}_h} (\estim{\sigma}^K)^2\Big)^\frac{1}{2}\\&=\estim{\mathrm{div}}+\estim{\mathrm{g}}+\estim{\sigma}
\end{align*}
 as the \emph{numerical component} of the estimator, and to $$\estim{\mathrm{def}}=\Big(\sum_{F \in \Fr} \alpha_3\estim{F}^2\Big)^\frac{1}{2}$$ as the \emph{defeaturing component}. The \emph{total estimator} is hence defined as
 $$\estim{}=\estim{\mathrm{num}}+\estim{\mathrm{def}}.$$
 \begin{rem}
 	{Let us remark that, compared to \cite{BCV2022_arxiv}, the use of an equilibrated flux reconstruction allows us to obtain a constant equal to 1 multiplying the bulk component of the numerical estimator $\estim{\sigma}$, and hence to bound sharply the numerical source of error far from included features. Note that  $\estim{\rm{num}}$ is not totally free of constants, but the terms involving $\alpha_1$ and $\alpha_2$ are localized uniquely on cut elements. Moreover, the equilibrated flux can be used also inside the defeaturing component of the estimator, instead of evaluating the normal trace of the numerical flux as proposed in \cite{BCV2022,BCV2022_arxiv}.}
 \end{rem}
\begin{rem}
Proposition \ref{prop1} states the reliability of the error estimator, but not its efficiency. Devising an efficient estimator would raise new challenges which are out of the scope of this paper, since we are interested in bounding the energy norm of the overall error in $\Omega$, while the numerical approximation error is committed in $\Omega^\star$, {with $\Omega \subset \Omega^\star$. Let us observe that, in presence of positive features only (see Remark \ref{rem:pos}), the defeatured domain is instead contained in $\Omega$, and the interface between the defeatured domain and a neglected positive feature plays a role similar to the artificial boundary introduced when truncating an unbounded domain. In that setting, the efficiency of an equilibrated flux estimator has been proven in \cite{ChaumontFrelet2025}, under the condition that the elements touching the artificial boundary are not too refined, and similar arguments could be used to address the efficiency of the defeaturing component of the estimator in the positive feature case. The extension to the unfitted case after feature reintroduction would however not be trivial and would require an analysis of its own.}
\end{rem}

\begin{rem}
Let us observe that the constants in Proposition \ref{prop1} are independent of how the feature boundaries intersect the mesh. Indeed, badly cut elements do not directly cause these constants to blow up; rather, their impact arises from the potential loss of stability of the flux reconstruction.
\end{rem}

	\section{The adaptive strategy} \label{sec:adaptive}
	In this section we aim at defining an adaptive refinement strategy based on the a posteriori error estimator proposed in Section \ref{sec:aposteriori}. As previously mentioned, we want to perform both numerical and geometrical adaptivity. Starting from a totally defeatured domain $\Omega_0$, from which a set of features $\mathcal{F}$ has been removed, and from a mesh $\mathcal{T}_h^0$, which is conforming to the boundary of $\Omega_0$ but completely blind to the features, we aim at pointing out:
	\begin{itemize}
		\item where and when the mesh needs to be refined ($h$-adaptivity);\vspace{0.3cm}
		\item which features need to be re-included in the geometry and when (geometrical adaptivity).
	\end{itemize}
 In the framework of adaptive finite elements for elliptic PDEs \cite{NochettoPRIMER}, at each iteration of the adaptive process we need to go though the following blocks:
\begin{equation*}
\boxed{\mathrm{SOLVE}}\rightarrow\boxed{\mathrm{ESTIMATE}}\rightarrow\boxed{\mathrm{MARK}}\rightarrow\boxed{\mathrm{REFINE}}
\end{equation*}
We now elaborate on each of these blocks, denoting by $s\in \mathbb{N}$ the current iteration index. 
We denote by $\Omega_{s}^\star$ the partially defeatured geometry at the beginning of the $s$-th iteration, choosing $\Omega_0^\star=\Omega_0$, and by $\mathcal{T}_h^{s}$ the mesh covering $\Omega_0$, obtained by the refinement of the starting mesh $\mathcal{T}_h^0$ through the $h$-adaptivity performed at the previous $s-1$ iterations. We remark again that the refined meshes are triangulations of $\Omega_0$, and therefore they are blind to the features.

	\subsection{Solve}\label{sec:solve}
	At step $s$, a problem in the form of \eqref{eq:prob_num} is solved, with $\Omegas=\Omega_s^\star=\text{int}\big(\overline{\Omega}\cup\lbrace\overline F,~F \in \tilde{\mathcal{F}}_s\rbrace\big) $, where $\tilde{\mathcal{F}}_s$ is the set of features still neglected (i.e. filled with material) at the beginning of iteration $s$, and which may be marked for geometric refinement at the end of the iteration. We define $\mathcal{F}_s^\star=\mathcal{F}\setminus\tilde{\mathcal{F}}_s$ as the set of features which have already been re-included into the geometry at the beginning of iteration $s$.  The active portion $\mathcal{A}_h^s$ of the mesh $\mathcal{T}_h^s$ is defined as
	$$\mathcal{A}_h^s=\lbrace K \in \mathcal{T}_h^s:~K \cap \Omega_s^\star \neq \emptyset\rbrace,$$  and $\Omega_\mathcal{T}=\Omega_\mathcal{T}^s=\bigcup_{K \in \mathcal{A}_h^s}\overline{K}$. The obtained discrete defeatured solution $\uhs=u_{h,s}^\star\in V_h(\OmegaT^s)$ is considered as an approximation of the solution $u$ of \eqref{eq:weak_omega} at the $s$-th iteration.
	
	\subsection{Estimate}
Let $\gr_0=\gr_{0,s}$, $\gr=\gr_s$, $\gamma_0^\star=\gamma_{0,s}^\star$ and $\gs=\gamma_s^\star$ be defined as in \eqref{gamma_star}-\eqref{gamma_tilde}, with $\tilde{\mathcal{F}}=\tilde{\mathcal{F}}_s$ and $\mathcal{F}^\star=\mathcal{F}_s^\star$ and let 
	$$\mathcal{G}_h^s=\lbrace K \in \mathcal{A}_h^s :~|K\cap\gamma_s^\star| \neq \emptyset\rbrace.
	$$
		Let $\flux=\bm{\sigma}_{h,s}$ be the (weakly) equilibrated flux computed from $u_{h,s}^\star$ by solving \eqref{eq:prob_eq_flux_nit} with $p=1$. For each $K \in \mathcal{A}_h^s$ we compute the quantities $\estim{\mathrm{div}}^K$, $\estim{\mathrm{g}}^K$ and $\estim{\sigma}^K$ defined in \eqref{estim_div_and_g}-\eqref{estim_sigma}.
	Let us introduce also
	$$\estim{K}^2:=\alpha_1(\estim{\mathrm{div}}^K)^2+\alpha_2(\estim{\mathrm{g}}^K)^2+(\estim{\sigma}^K)^2,$$
	and let us observe that $\estim{K}^2=(\estim{\sigma}^K)^2$ if $K \in \mathcal{A}_h^s \setminus \mathcal{G}_h^s$.
Finally, for each feature $F\in \Fr_s$,
 we compute $\estim{F}$ as defined in \eqref{eq:estim_gamma}. 

\subsection{Mark}\label{subsec:mark}	
In order to mark the elements in $\mathcal{A}_h^s$ and the features in $\Fr_s$ we adopt a D\"{o}rfler strategy \cite{Dorfler}.  Let us fix a marking parameter $\theta\in (0,1]$ and let $\mathcal{M}^s=\mathcal{A}_h^s \cup \Fr_s$. For each element or feature $M\in \mathcal{M}^s$ we define
		$$\mathcal{E}_M^2=\begin{cases} \estim{K}^2 & \text {if } M=K\in\mathcal{A}_h^s\\
		\alpha_3\estim{F}^2 &\text {if } M =F\in \Fr_s. \end{cases}$$ We aim at marking the smallest subset $\mathcal{D}^s \subset \mathcal{M}^s$ such that
		\begin{equation}
		\sum_{M \in \mathcal{D}^s} \mathcal{E}_M^2\geq\theta \sum_{M\in \mathcal{M}^s}\mathcal{E}_M^2. \label{Dorfler}
		\end{equation}
		\subsection{Refine}
\begin{figure}
					\centering
					\begin{subfigure}[t]{.32\textwidth}
						\centering
	\includegraphics[width=1\linewidth]{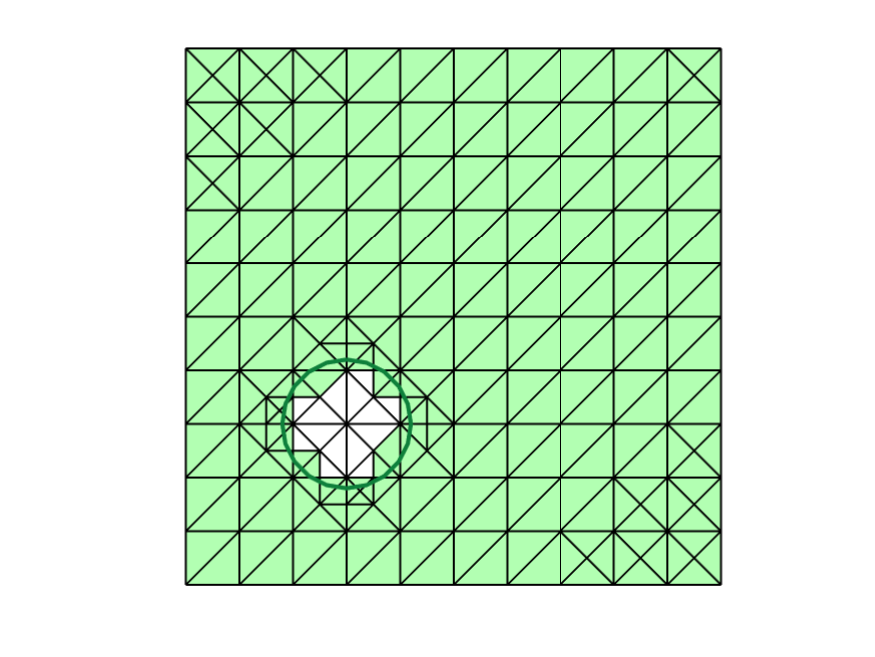}
		\caption{Active mesh $\mathcal{A}_h^s$}
	\label{fig:adaptivity1}
					\end{subfigure}
				\begin{subfigure}[t]{.32\textwidth}
				\centering
				\includegraphics[width=1\linewidth]{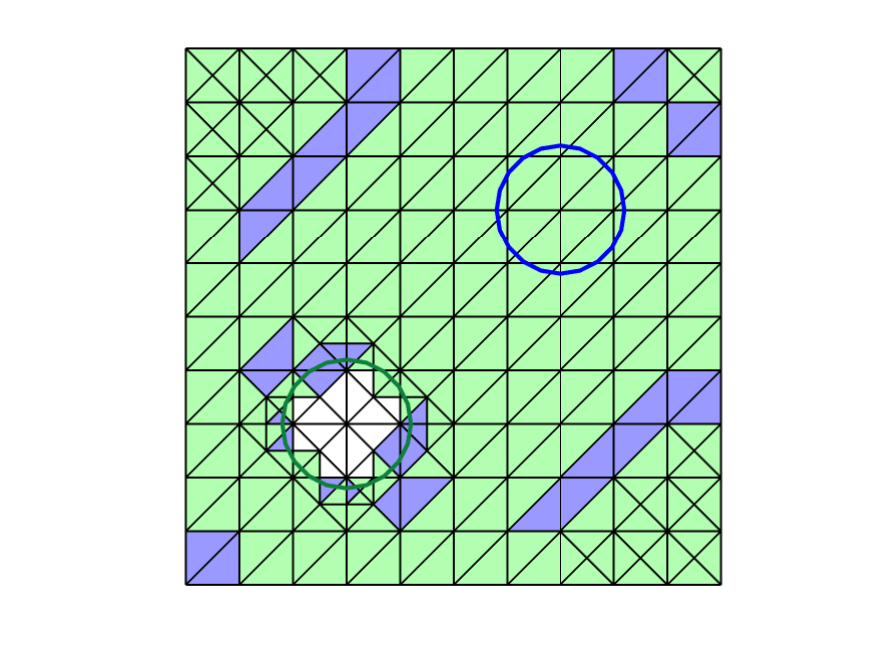}
				\caption{Feature and elements marked for refinement (in blue)}
				\label{fig:adaptivity2}
			\end{subfigure}
			\begin{subfigure}[t]{.32\textwidth}
			\centering
			\includegraphics[width=1\linewidth]{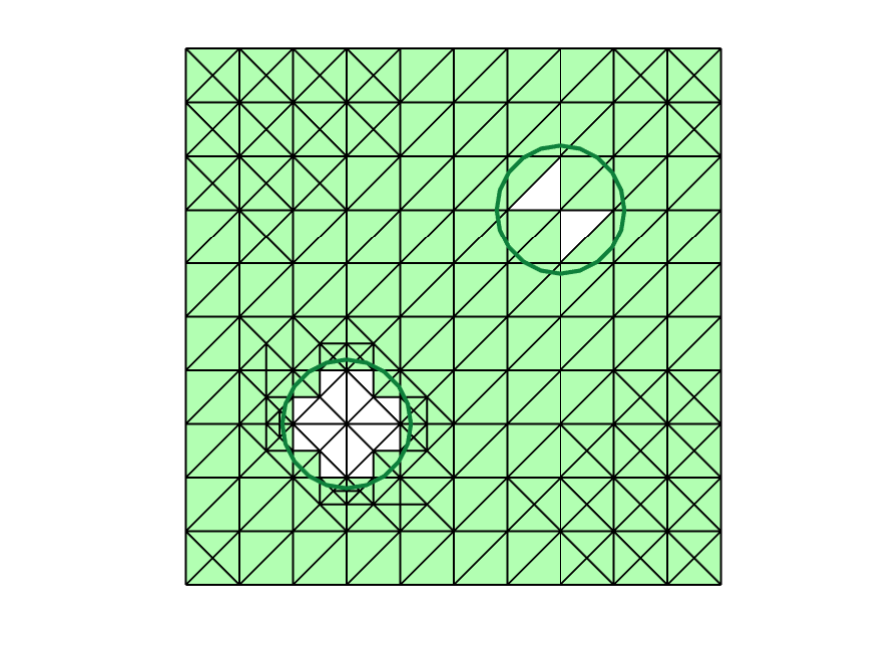}
			\caption{Active mesh $\mathcal{A}_h^{s+1}$}
			\label{fig:adaptivity3}
		\end{subfigure}
	\caption{Example of an iteration of combined adaptivity}
	\label{fig:adaptive}
				\end{figure}

		The mesh is refined thanks to an $h$-refinement procedure applied to the elements marked in $\mathcal{D}^s \cap \mathcal{A}_h^s$, leading to a mesh $\mathcal{T}_h^{s+1}$. The marked features in $\mathcal{D}^s \cap\Fr_s$ are instead included into the geometry, i.e.
				$$\Omega_{s+1}^\star=\Omega_{s}^\star\setminus\bigcup_{F \in \mathcal{D}^s\cap \Fr_s}\overline{F}$$ and $\Fr_{s+1}=\Fr_s\setminus (\mathcal{D}^s \cap \Fr_s)$, $ \mathcal{F}_{s+1}^\star=\mathcal{F}_{s}^\star\cup (\mathcal{D}^s \cap \Fr_s) .$
				The active mesh for the following step is defined as $\mathcal{A}_h^{s+1}=\lbrace K \in \mathcal{T}_h^{s+1}:~K\cap \Omega_{s+1}^\star \neq \emptyset\rbrace$ (see Figure~\ref{fig:adaptive}).
			\begin{rem}\label{rem:assum}
				It is interesting to remark how, similarly to \cite{BCV2022_arxiv}, the proposed adaptive strategy naturally takes care of Assumption \ref{assum:hf} {and also of Assumption~4}.
				
				 Let $F\in\Fr_s$, i.e. $F$ is a feature that is not included into the geometry at step $s$. If $F\cap K=\emptyset$ $\forall K \in \mathcal{G}_h^s$, which means that $F$ is not intersecting any element which is also intersected by a feature in $\Fs_s$, then Assumption \ref{assum:hf} is actually unnecessary for the given feature, since $||\nabla \cdot \bm{\sigma}_{h,s}-f||_F=0$ (see the proof of Proposition \ref{prop1}). However, it will still be naturally satisfied since, if the error is concentrated in the feature, this will lead to the inclusion of $F$ into the geometry.
				 The same holds if there exists at least one element $K \in \mathcal{G}_h^s$ such that $K\cap F\neq \emptyset$. Indeed, the local contribution to the numerical component of the estimator in the elements in $\mathcal{G}_h^s$ is higher, due to the terms $\estim{\mathrm{div}}^K$ and $\estim{\mathrm{g}}^K$. This means that $K$ will be early marked for refinement, allowing to separate $F$ from the elements in $\mathcal{G}_{h}^s$ in a few iterations. 
				 
				 {If we consider, instead, $F \in \mathcal{F}_s^\star$, i.e. $F$ is a feature which has already been re-included into the geometry at the beginning of iteration $s$, Assumption \ref{assum4} tells us that the size of the elements intersecting the boundary of $F$ should be controlled from above by the size of the feature itself. This is handled naturally by the adaptive strategy as well. Indeed, once the feature is marked for re-inclusion, the terms $\estim{\mathrm{div}}^K$ and $\estim{\mathrm{g}}^K$ increase the local contribution to the numerical component of the estimator in the elements cut by the boundary of $F$, hence leading to their refinement and to the fulfillment of the condition in a few iterations.}
				\end{rem}

\section{Numerical experiments}\label{sec:num_exp}
In the following we propose some numerical examples in order to validate the proposed estimator. We focus on the $r=1$ case, i.e. we use linear Lagrange finite elements to compute the primal solution $\uhs$, and we reconstruct the (weakly) equilibrated flux in a Raviart--Thomas space of order 1, {i.e. with\vspace{-0.2cm} $$\bm{M}_h(\wa):=\left\lbrace  \bm{v}_h\in \Hdiv{\wa}:~\tr{\bm{v}_h}{K} \in [\mathcal{P}_1(K)]^d+\bm{x}\mathcal{P}_1(K),~\forall K \in \wa\right\rbrace\vspace{-0.2cm}$$ in \eqref{eq:Mha}}. Some remarks on the stabilized version (see Problem \eqref{eq:prob_eq_flux_nit_stab}) are reported for the first numerical example. For the elements $K \in \mathcal{G}_h^s$, the integrals on the active portion $K\cap \Omega_s^\star$ are computed by defining proper quadrature rules on a subtriangulation of the active portion of the element itself. The integrals on $\gamma_{\bm{a}}^\star$ appearing in \eqref{eq:prob_eq_flux_nit} are computed by quadrature formulas on the intersection between $\mathcal{A}_h$ and the feature edges (see for instance \cite{Boffi_nonmatching} for a review on integration of interface terms in non-matching techniques). The same is done also in the computation of the integrals on $\tilde{\gamma}$ required in the evaluation of the defeaturing component of the estimator.

Let us remark that when no feature is included into the geometry, we have $\estim{K}^2=\estim{\sigma}^2$ for all $K \in \mathcal{A}_h^s$ and for all $s \in \mathbb{N}$. Indeed, $\estim{\mathrm{g}}=0$ and, since $f\in Q_h(\Omega_0)$, also $\estim{\mathrm{div}}=0$. Hence, when in the following numerical examples we will compare the performances of the proposed adaptive procedure with standard mesh refinement, we will refer to a standard $h$-adaptivity in which the mesh is refined by the D\"orfler strategy reported in Section \ref{subsec:mark} according only to the value of $(\estim{\sigma}^K)^2=||\flux+\nabla u_h^0||^2_{K}$, $K\in \mathcal{T}_h^s$. We will refer to the adaptive procedure allowing both for mesh and geometric refinement as \emph{combined} adaptivity. In the tests that follow, convergence trends are shown with respect to the number $N$ of degrees of freedom of $\uhs$. The numerical component of the estimator is approximated as $\estim{\mathrm{num}}=\big(\sum_{K \in \mathcal{A}_h}\estim{K}^2\big)^\frac{1}{2}$.

\subsection{Test 1: single internal feature}
\begin{figure}
 	\centering
 	\begin{subfigure}{.5\textwidth}
 		\centering
 		\includegraphics[width=0.74\linewidth]{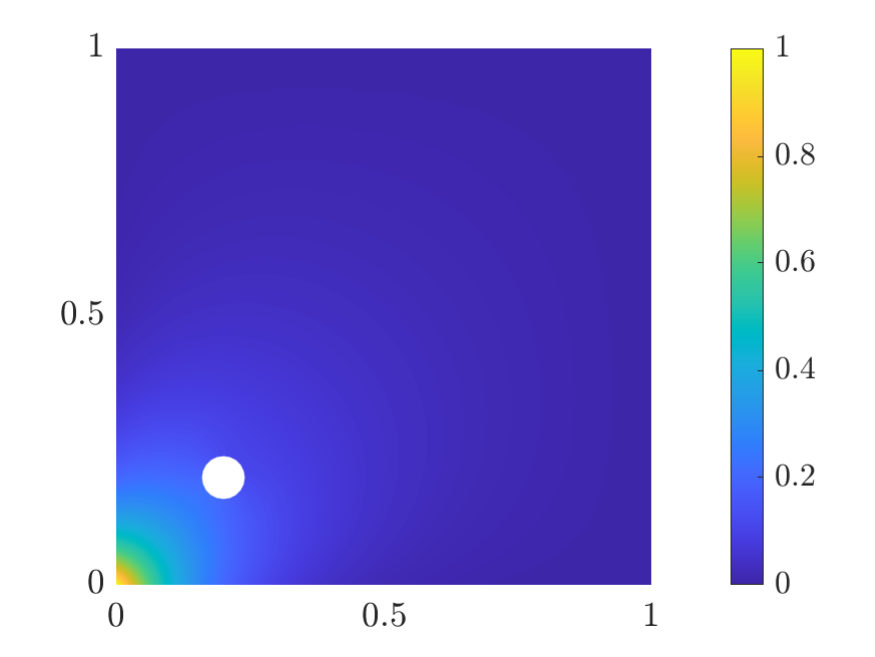}
 		\caption{Reference solution}
 		\label{fig:t1_meshsol-a}
 	\end{subfigure}%
 	\begin{subfigure}{.5\textwidth}
 		\centering
 		\includegraphics[width=0.74\linewidth]{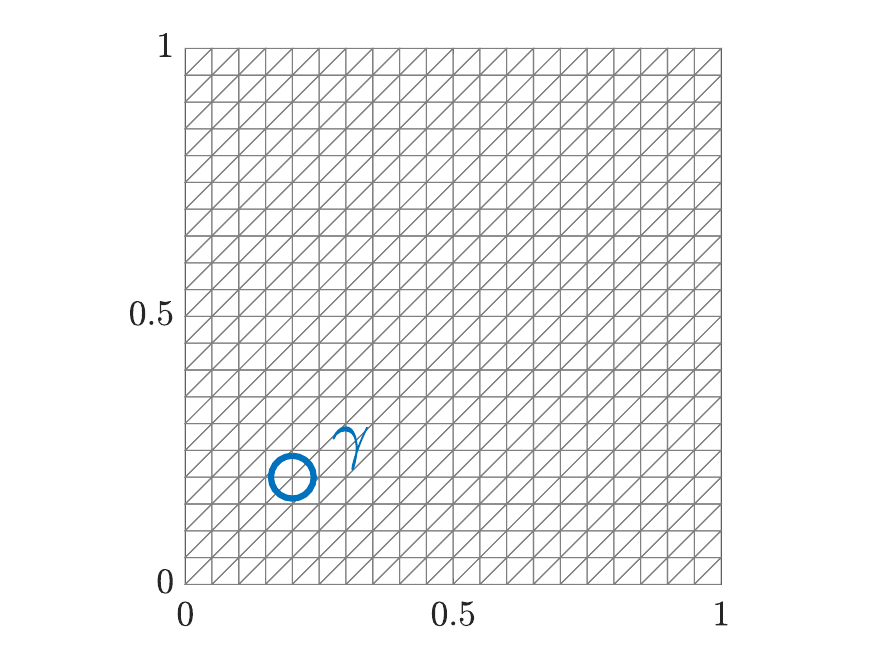}
 		\caption{Non conforming initial mesh on $\Omega_0$}
 		\label{fig:t1_meshsol-b}
 	\end{subfigure}
 	\caption{Test1: reference solution and initial mesh on $\Omega_0$.}
 	\label{fig:t1_meshsol}
 \end{figure}

Let $\Omega_0=(0,1)^2$ and let $F$ be a regular polygonal feature with 20 faces, having boundary $\gamma$, inscribed in a circle of radius $\varepsilon=0.04$ and centered in $(0.2,0.2)$ (see Figure~\ref{fig:t1_meshsol}). 
Let $\Omega=\Omega_0\setminus \overline{F}$ and let $\Gamma_D=\partial \Omega_0$. We consider the following problem:
\begin{equation}
\begin{cases}
-\Delta u=0 &\text{in } \Omega\\
u=g_{\mathrm{D}} &\text{on } \Gamma_{\mathrm{D}}\\
\nabla u\cdot \bm{n}=0 &\text{on } \gamma
\end{cases}\label{prob_test}
\end{equation}
with 
\begin{equation}
g_{\mathrm{D}}(x,y)=\begin{cases} e^{-8(x+y)} &\text{if } x=0 \text{ or } y=0\\
0 &\text{otherwise}
\end{cases} \label{test:dirichlet}
\end{equation}
 The aim of this first numerical example is to analyze the convergence of the combined adaptive strategy, comparing it to the convergence of a standard one, which performs only mesh refinement. The error is computed with respect to a reference solution reported in Figure \ref{fig:t1_meshsol-a}, obtained on a very fine mesh conforming to the boundary of the feature. For what concerns the total estimator $\estim{}$, we choose $\alpha_1=\alpha_2=\alpha_3=1$. Both strategies start from the mesh reported in Figure \ref{fig:t1_meshsol-b}, ($h=7.07\cdot 10^{-2}$, $N=361$ degrees of freedom for the variable $\uhs$), and in both cases we choose $\theta=0.3$ in \eqref{Dorfler}. The adaptive strategy is stopped as $5\cdot 10^3$ degrees of freedom are reached for $\uhs$.
 \begin{figure}
	\centering
		\includegraphics[width=0.48\linewidth]{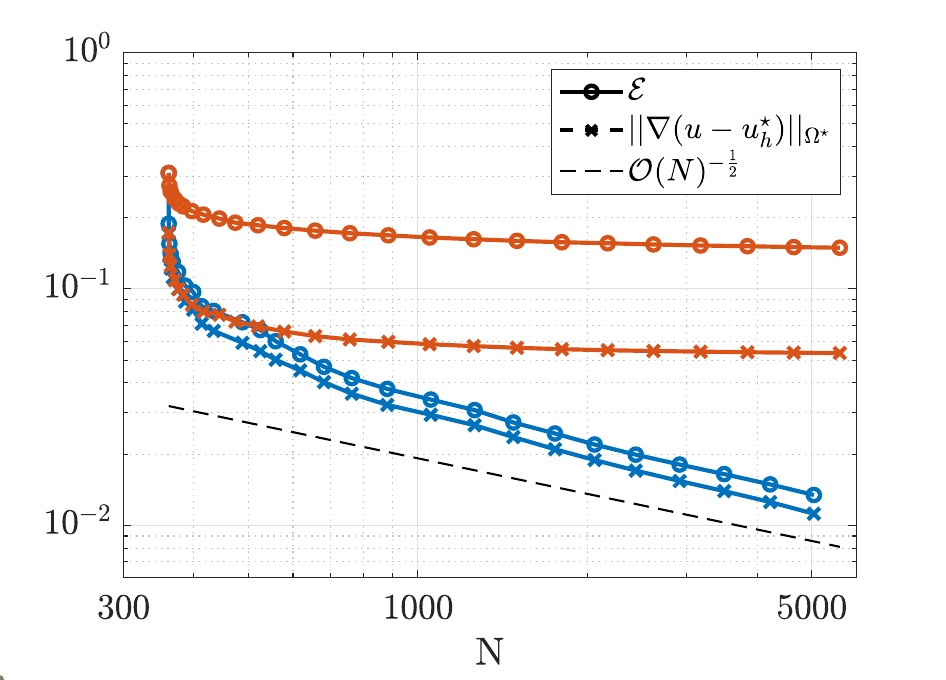}\vspace{-0.5cm}
	\caption{Test1: Convergence of the total estimator and of the overall error against the number of degrees of freedom ($N$). In blue: feature included at the first iteration; in red: feature never included.}
	\label{fig:t1:convergence}
\end{figure}
 
 The trend of the overall error and of the total error estimator during the adaptive procedure is reported in Figure \ref{fig:t1:convergence}. 
 The red curves refer to the case in which the feature is not included into the geometry and only mesh refinement is performed, whereas the blue curves are obtained allowing both for mesh and geometric refinement. We can observe that, if the feature is not included, the error and the total estimator reach a plateau. On the contrary, if the feature is included into the geometry the error and the total estimator converge at the expected rate, i.e. as $\mathcal{O}(N^{-\frac{1}{2}})$, where $N$ is the total number of degrees of freedom for the variable $\uhs$. In particular, following the procedure devised in Section \ref{subsec:mark}, the feature is included at the very first step, due to the predominance of the defeaturing component of the estimator on the local contributions to the numerical component. Hence, excluding the very first iteration, the whole adaptive procedure is steered by the numerical component of the estimator, as in the case in which the feature is not included. However, when the feature is included, the numerical component of the estimator accounts also for the error introduced by the weak imposition of the Neumann condition on the feature boundary. 
 The feature is included at the very first iteration since it is located in a region in which the gradient of the solution is very steep: hence it is expected to have a strong impact on the accuracy of the solution. The step at which the feature is included could be slightly postponed by decreasing the value of $\alpha_3$, within the limit in which the defeaturing component of the estimator does not affect the convergence of the total estimator. As discussed in \cite{BCV2022}, the optimal value of $\alpha_3$ results to be problem dependent and the parameter is usually fine tuned heuristically.
 
\begin{figure}
	\centering
	\begin{subfigure}{.5\textwidth}
		\centering
		\includegraphics[width=0.83\linewidth]{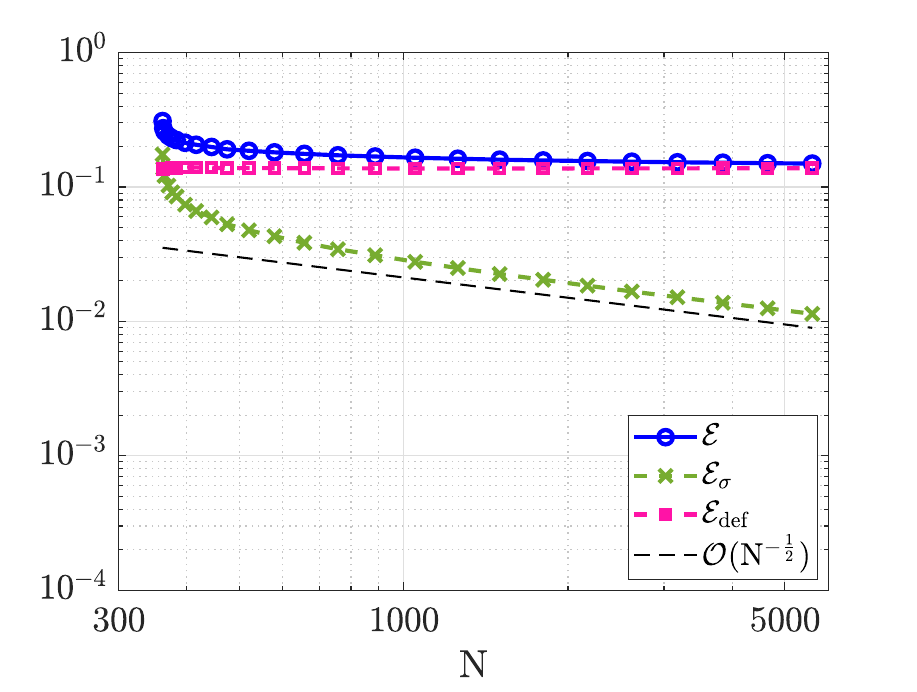}
		\caption{No feature is included}
		\label{fig:t1:converegnceNOADD}
	\end{subfigure}%
	\begin{subfigure}{.5\textwidth}
		\centering
		\includegraphics[width=0.83\linewidth]{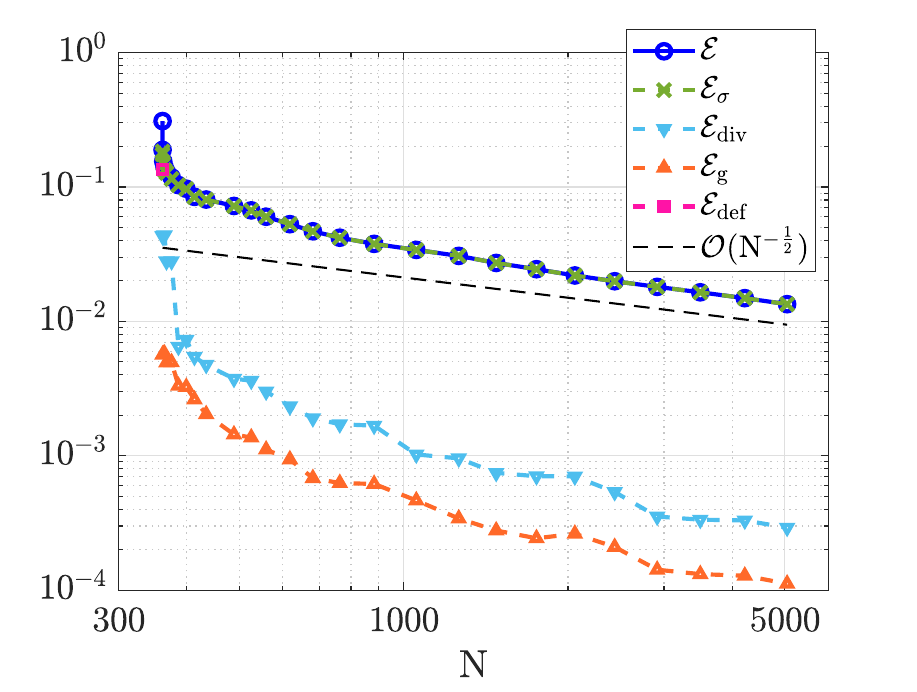}
		\caption{Feature included at iteration 1}
		\label{fig:t1:convergenceADD1}
	\end{subfigure}
	\caption{Test1: Convergence of the total estimator and of its components with respect to the number of degrees of freedom ($N$).}
	\label{fig:t1:convergence_detail}
\end{figure}

The detail on the convergence of the components of the total estimator is reported in Figure~\ref{fig:t1:convergence_detail}. As expected, when the feature is not added (Figure~\ref{fig:t1:converegnceNOADD}) the defeaturing component $\estim{\mathrm{def}}$ of the estimator remains almost constant, and for this reason, even if $\estim{\mathrm{num}}=\estim{\sigma}$ converges at the expected rate (as $N^{-\frac{1}{2}}$), the total estimator reaches a plateau. On the contrary, when the feature is added, $\estim{\mathrm{def}}=0$ starting from the second iteration (Figure~\ref{fig:t1:convergenceADD1}) and the total estimator converges as $N^{-\frac{1}{2}}$.

\begin{figure}
	\centering
	\begin{subfigure}{.5\textwidth}
		\centering
		\includegraphics[width=0.77\linewidth]{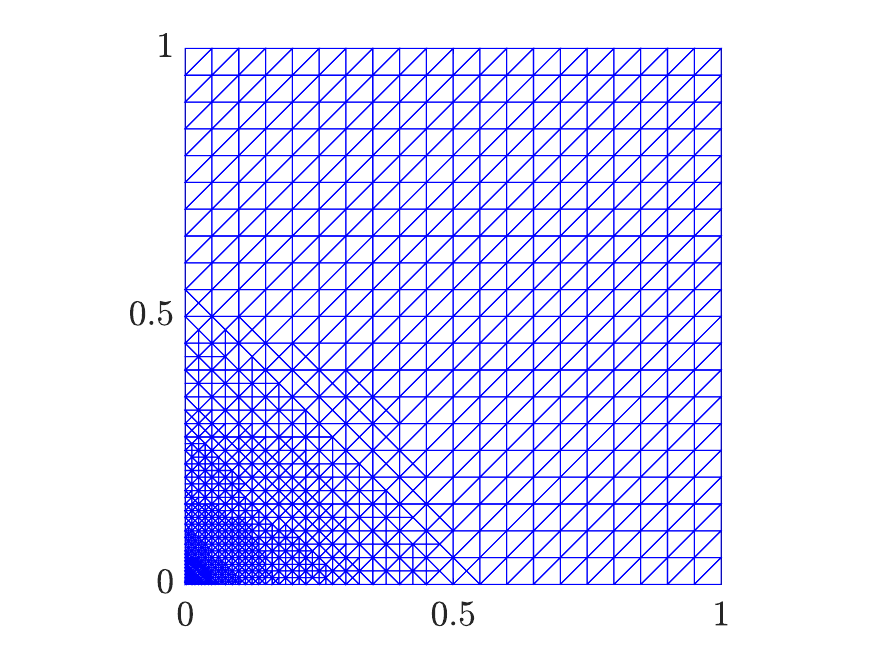}
		\caption{No feature is included}
		\label{fig:mesh_noAdd}
	\end{subfigure}%
	\begin{subfigure}{.5\textwidth}
		\centering
		\includegraphics[width=0.77\linewidth]{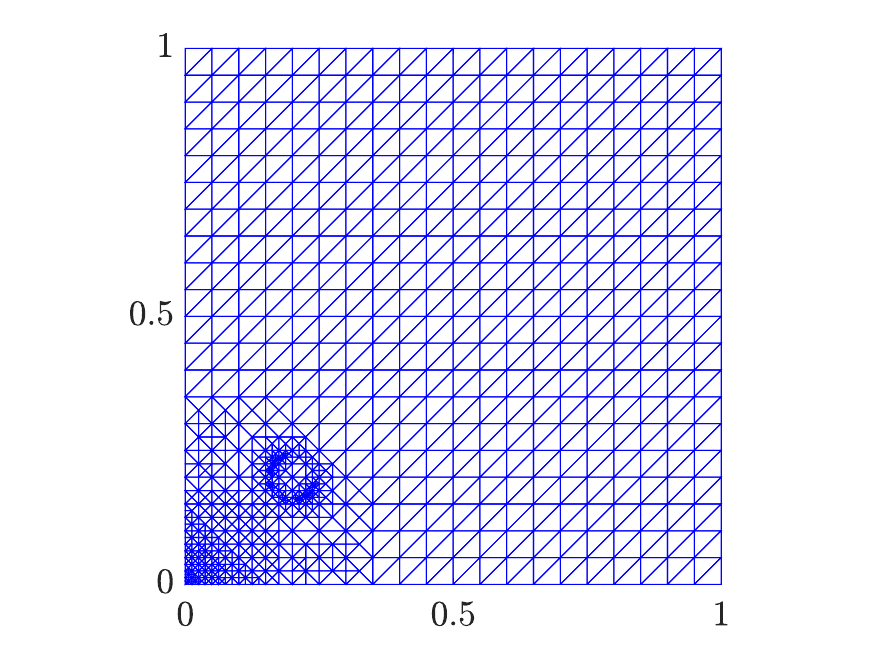}
		\caption{Feature is included at iteration 1}
		\label{fig:mesh_add1_noP1}
	\end{subfigure}
	\caption{Test 1: Adapted mesh after 15 steps of the adaptive procedure.}
	\label{fig:t1:mesh16}
\end{figure}
 
 The meshes obtained after 15 steps of the adaptive procedure are reported in Figure \ref{fig:t1:mesh16}, both for standard mesh adaptivity and for combined adaptivity. As expected, in both cases the mesh appears to be refined towards the bottom left corner of the domain. However, we can observe how the combination of mesh and geometric refinement allows the mesh to capture the feature boundary.
 
\begin{figure}
	\centering
	\begin{subfigure}{.5\textwidth}
		\centering
		\includegraphics[width=0.8\linewidth]{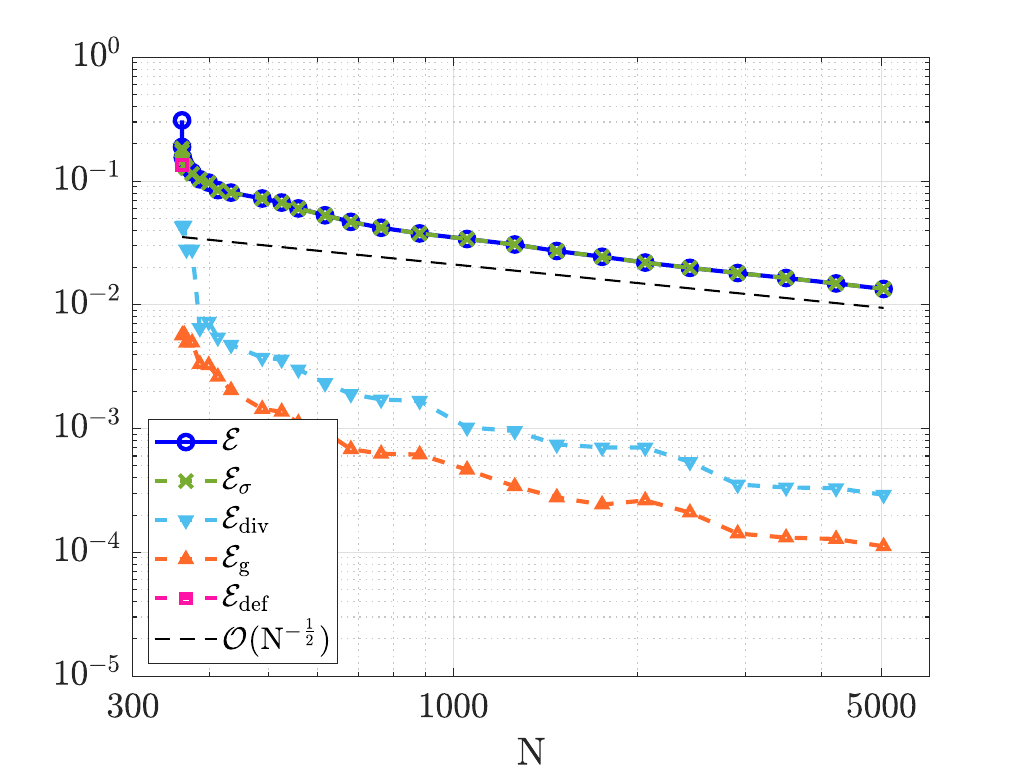}
		\caption{No stabilization}
		\label{fig:conv_stab_a}
	\end{subfigure}%
	\begin{subfigure}{.5\textwidth}
		\centering
		\includegraphics[width=0.8\linewidth]{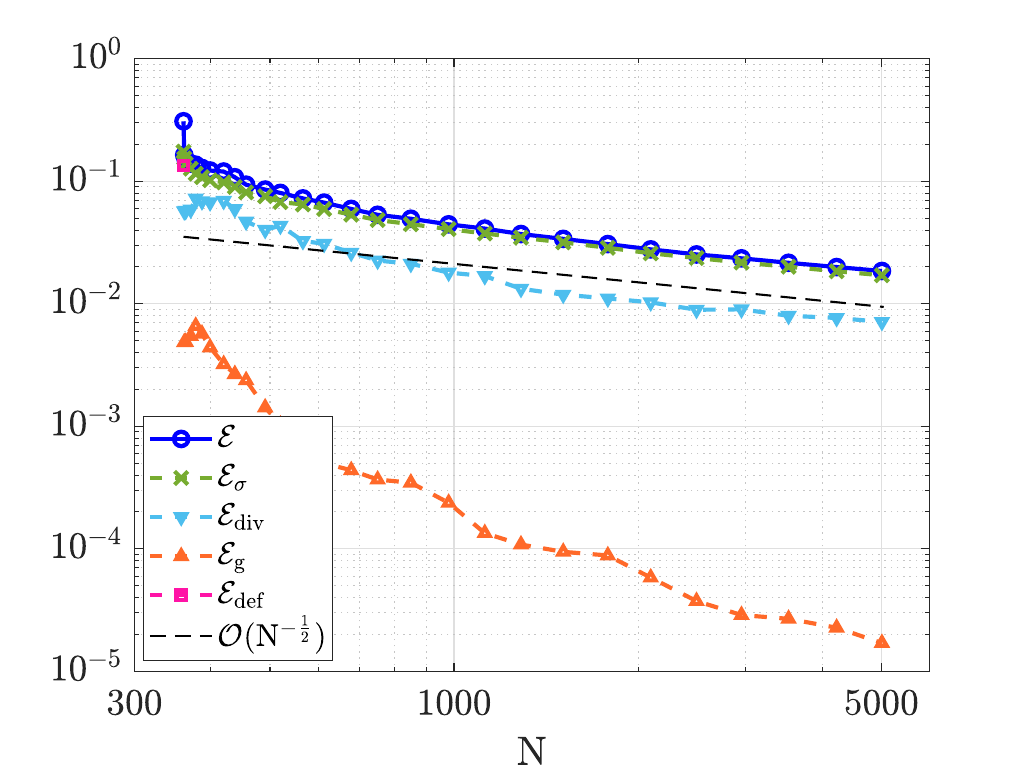}
		\caption{With stabilization}
		\label{fig:conv_stab_b}
	\end{subfigure}
	\caption{Test 1: Convergence with respect to the number of degrees of freedom (N) of the total estimator and of its components with and without stabilization}
	\label{fig:conv_stab}
\end{figure}

 As reported in Remark \ref{stab}, the weak imposition of Neumann boundary conditions in the computation of the weakly equilibrated flux reconstruction $\flux$ may require stabilization in presence of very small cuts. Figure \ref{fig:conv_stab_b} reports the detail on the convergence of the components of the estimator when the stabilized problem \eqref{eq:prob_eq_flux_nit_stab} is solved instead of \eqref{eq:prob_eq_flux_nit}. The curves obtained in the case with no stabilization, which were reported in Figure \ref{fig:t1:convergence_detail}, are re-proposed and rescaled in Figure \ref{fig:conv_stab_a} to ease the comparison.  As expected, the main drawback of activating the stabilization terms \eqref{stab1}-\eqref{stab2} is represented by an increase of $\estim{\mathrm{div}}$, which is due to the fact that mass conservation is further deteriorated with respect to the non-stabilized case, involving all the cut patches. However, this local increase of $\estim{\mathrm{div}}^K$ produces a stronger mesh refinement close to the feature, which ends up in a faster decrease of $\estim{\mathrm{g}}$. For this reason, the final impact of stabilization on the total estimator is pretty low, and by decreasing the value of $\alpha_{1}$ the effectivity index of the non stabilized case could be restored. Let us remark that in this case, in which the feature is added at the very first iteration, the inclusion of stabilizing terms does not deteriorate the convergence rate of the total estimator. In a different scenario, characterized by the presence of multiple features which are added at different stages of the adaptive procedure, the parameters $\alpha_{1},\alpha_{2},\alpha_3$ and $\beta_1,\beta_2$ would need to be properly tuned, which is not always trivial. For this reason, from a barely practical standpoint, the easiest way to avoid the issues related to small cuts, may be to neglect badly cut patches, which have a small impact on the reconstructed flux. Let us recall, indeed, that the reconstructed flux is needed only to evaluate the estimator, and is not the primal variable of the problem.
 
\subsection{Test 2: multiple internal and boundary features}
The aim of this second numerical experiment is to test the proposed adaptive procedure on a case characterized by the presence of multiple features, both internal and on the boundary. Let again $\Omega_0=(0,1)^2$ and let us consider a set $\mathcal{F}=\lbrace F_i\rbrace_{i=1}^{37}$ of 37 features, being all regular polygons inscribed in a circle of radius $\epsilon_i$, centered in $\bm{x}_i^C$, having $n_i^e$ edges and rotated of an angle $\theta_i$ with respect to a reference orientation in which one of the vertexes is located in $\bm{x}_i^C+[0,\epsilon_i]$. In particular, denoting by $\mathcal{U}(a,b)$ a standard uniform distribution between $a$ and $b$, we take $\epsilon_i\sim \mathcal{U}(1e-3,5e-2)$ and $\theta_i\sim \mathcal{U}(0,2\pi)$. For internal features $\bm{x}_i^C\sim\mathcal{U}(0.1,0.9)^2$, while for boundary features $\bm{x}_i^C\sim\lbrace0,1\rbrace \times\mathcal{U}(0.1,0.9)$. For both cases the number of edges $n_i^e$ is chosen from a discrete integer uniform distribution between 4 and 16. The actual data used in the experiments are available in the Appendix.

On $\Omega=\Omega_0\setminus \lbrace \overline F,~F \in \mathcal{F}\rbrace$ we consider the following test problem:
\begin{equation}
\begin{cases}
-\Delta u=0 &\text{in } \Omega\\
u=g_{\mathrm{D}} &\text{on } \Gamma_{\mathrm{D}}\\
\nabla u \cdot \bm{n}=0&\text{on } \Gamma_{\mathrm{N}}
\end{cases}\vspace{-0.2cm}
\end{equation}
with $\Gamma_{\mathrm{D}}=\lbrace(x,y):y\in \lbrace0,1 \rbrace\rbrace$ and\vspace{-0.2cm}
\begin{equation}
g_{\mathrm{D}}(x,y)=\begin{cases} e^{-8(x+y)} &\text{if } y=0 \\
0 &\text{if } y=1.
\end{cases}
\end{equation}
The behavior of $u$ is very similar to the one in Test 1; however the presence of Neumann boundaries allows us to set some features along the boundary itself. 
For what concerns the (partially) defeatured problem, we choose $g_0=0$ (see \eqref{eq:strong_omega0} and \eqref{eq:strong_omega_star}).
\begin{figure}
	\centering
	\includegraphics[width=0.45\linewidth]{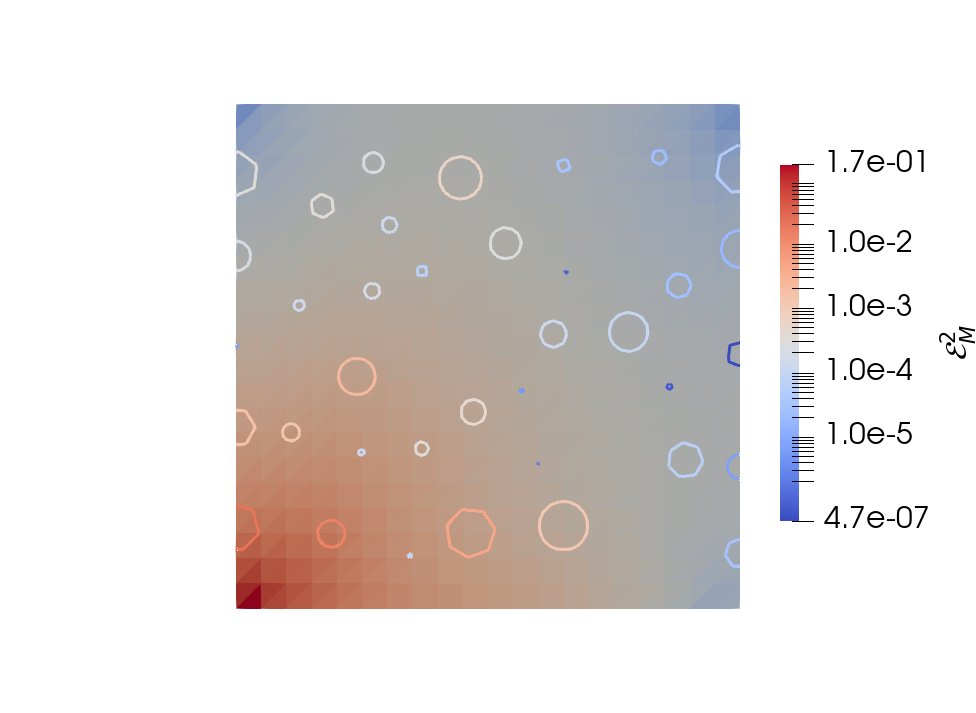}\vspace{-0.5cm}
	\caption{Test 2: Square of the local numerical components and of the local defeaturing components of the estimator at the first iteration of the adaptive procedure.}
	\label{fig:estimmstep1}
\end{figure}
The features in $\mathcal{F}$ are shown in Figure \ref{fig:estimmstep1}, which also reports the initial values of the square of the local contributions to the numerical and the defeaturing component of the estimators, i.e. the values contained in $\estim{M}^2$ at the very first application of the D\"orfler strategy reported in Section \ref{subsec:mark}. The initial mesh $\mathcal{T}_h^0$ is the same used in Test 1 ($h=7.07\cdot 10^{-2}$), ending however up in $N=399$ degrees of freedom for the variable $\uhs$, due to the different choice of the boundary conditions.

\begin{figure}
	\centering
	\includegraphics[width=0.5\linewidth]{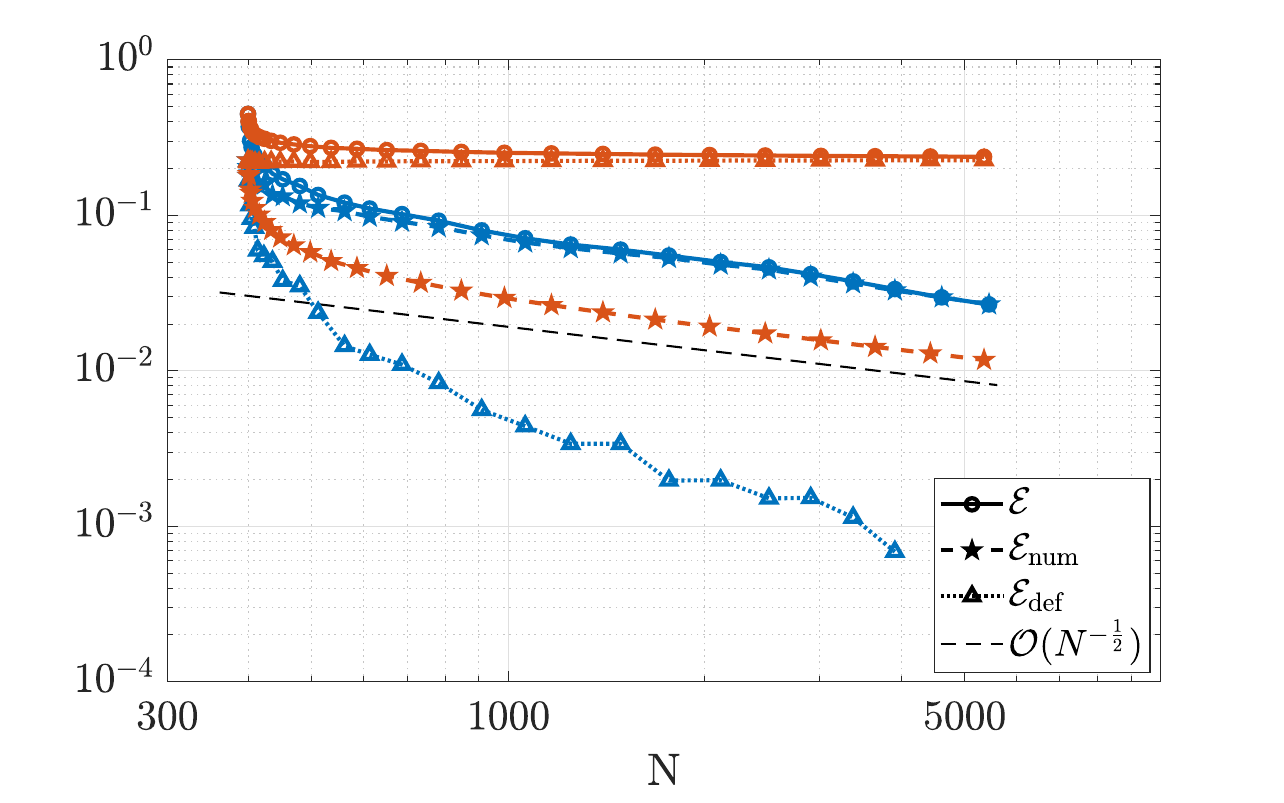}
	\caption{Test 2: Convergence of the total estimator and of its components at the increase of the number of degrees of freedom (N). In red, no feature included; in blue features included with $\alpha_3=1$.}
	\label{fig:t2:convergence}
\end{figure}

\begin{figure}
	\centering
	\begin{subfigure}{.5\textwidth}
		\centering
		\includegraphics[width=0.8\linewidth]{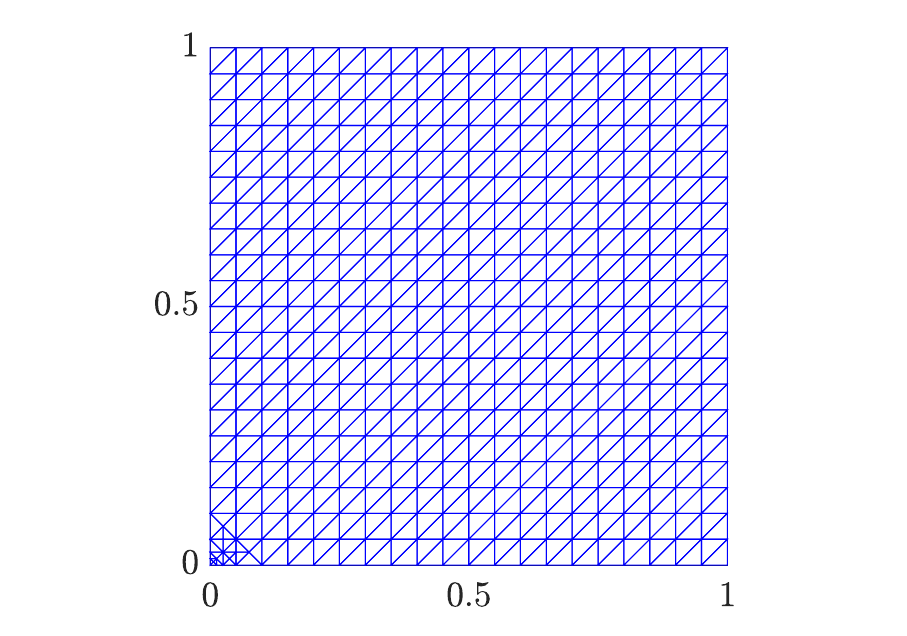}
		\caption{Mesh adaptivity}
		\label{fig:t2:mesh_after6it_noadd}
	\end{subfigure}%
	\begin{subfigure}{.5\textwidth}
		\centering
		\includegraphics[width=0.8\linewidth]{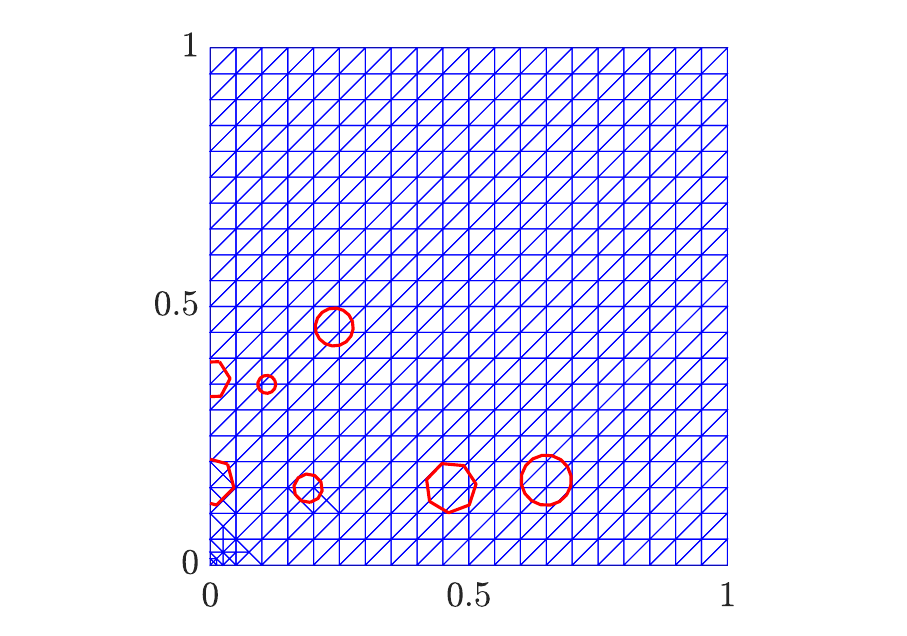}
		\caption{Combined adaptivity}
		\label{fig:t2:mesh_after6it_add}
	\end{subfigure}
	\caption{Test 2: Adapted mesh after 6 iterations of the adaptive procedure (included features highlighted for the geometric adaptivity case).}
	\label{fig:t2:mesh_after6it}
\end{figure}

\begin{figure}
	\centering
	\begin{subfigure}{.5\textwidth}
		\centering
		\includegraphics[width=0.8\linewidth]{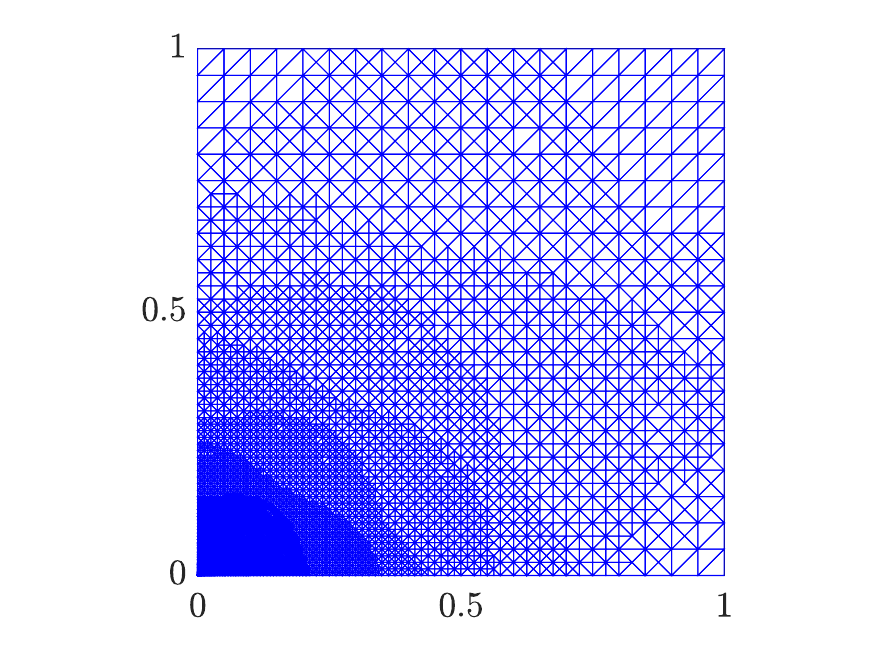}
		\caption{Mesh adaptivity}
		\label{fig:t2:mesh_final_add}
	\end{subfigure}%
	\begin{subfigure}{.5\textwidth}
		\centering
		\includegraphics[width=0.8\linewidth]{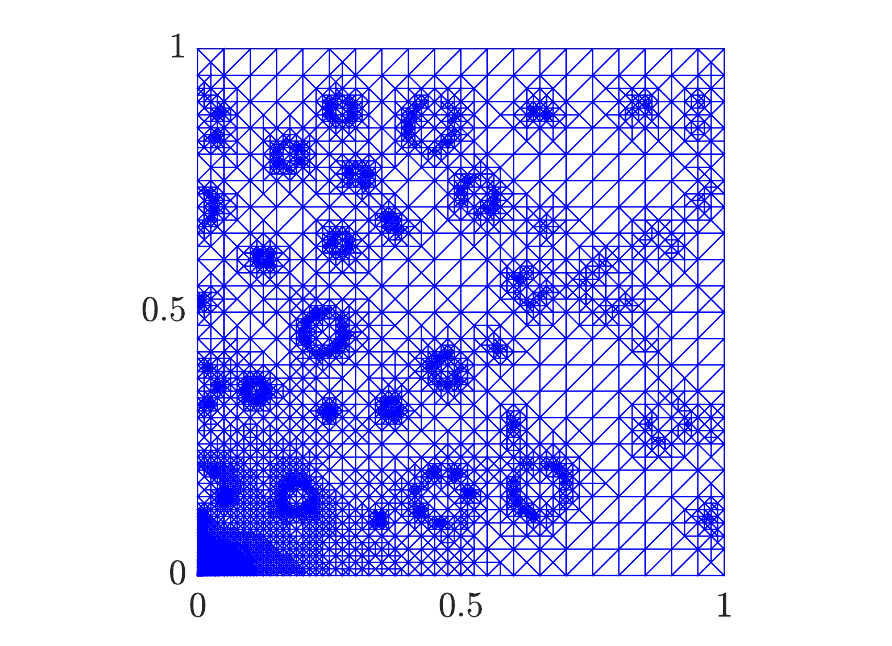}
			\caption{Combined adaptivity}
		\label{fig:t2:mesh_final_noadd}
	\end{subfigure}
	\caption{Test 2: Adapted mesh after 27 iterations of the adaptive procedure.}
	\label{fig:t2:mesh_final}
\end{figure}

Figure \ref{fig:t2:convergence} compares the convergence of the error estimator and of its components in two scenarios: the blue curves refer to the case in which the adaptive strategy proposed in Section \ref{sec:adaptive} is applied with $\alpha_1=\alpha_2=\alpha_3=1$, whereas the red curves are obtained without including any feature, i.e. with standard mesh refinement based on the value of $(\estim{\sigma}^K)^2$. As in the previous numerical test, we can observe how, in this second case, the total estimator rapidly reaches a plateau, due to the fact that the defeaturing component of the estimator remains constant. Allowing instead for both mesh and geometric refinement the total estimator tends to converge as $N^{-\frac{1}{2}}$, i.e. as the numerical component. Figure \ref{fig:t2:mesh_after6it} shows the computational mesh after 6 iterations of the adaptive procedure, both for standard mesh adaptivity and for the combined one. At this stage, if no features are included, the total estimator is reduced of about $27\%$ with respect to its initial value while, if we allow also for geometric adaptivity, $\estim{}$ is reduced of about $50\%$, with only 7 features included. As expected, and as shown in Figure \ref{fig:t2:mesh_after6it_add}, these 7 features are the biggest among the ones closest to the bottom left corner, which is coherent also with the results in \cite{BCV2022_arxiv}. The adaptive process is stopped as $5\cdot 10^3$ degrees of freedom are reached. At this stage, in the combined strategy, all the features are included (all the features are actually included before iteration 26) and the total estimator has been reduced of about $94\%$, against a relative decrease of $47\%$ in the standard mesh adaptivity case. Figure \ref{fig:t2:mesh_final} shows the final meshes, both for standard and combined adaptivity.


\subsection{Test 3: multiple internal features, discontinuous coefficient}
\begin{figure}
	\centering
		\begin{subfigure}{.5\textwidth}
			\centering
		\includegraphics[width=0.8\linewidth]{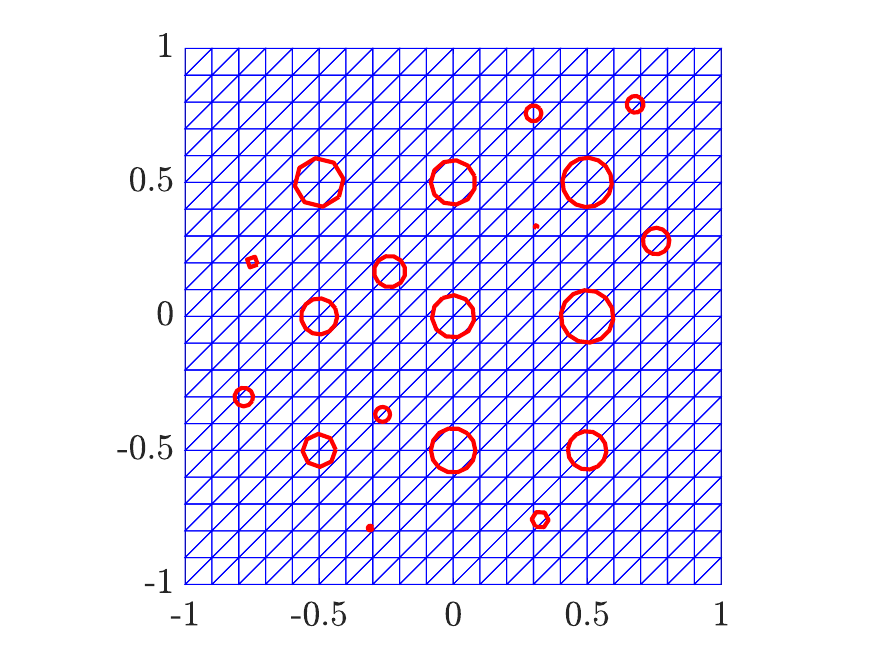}
		\caption{Features and initial mesh}
		\label{fig:t3:initialMesh}
	\end{subfigure}%
	\begin{subfigure}{.5\textwidth}
		\centering
		\includegraphics[width=0.8\linewidth]{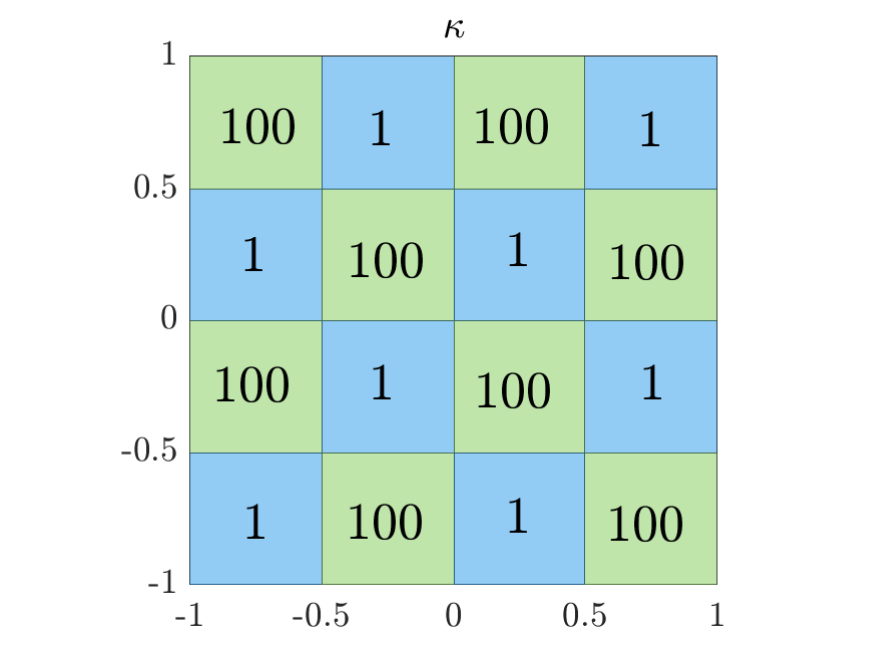}
	\caption{Diffusion coefficient $\kappa$}
	\label{fig:diffusiongrid}
	\end{subfigure}
	\caption{Test 3: initial mesh, features and chessboard diffusion coefficient}
\label{fig:t3:mesh_and_diff}

\end{figure}

For this last numerical example let us choose $\Omega_0=(-1,1)^2$ and let $\mathcal{F}$ be a set of 19 negative, internal, polygonal features (see Figure \ref{fig:t3:initialMesh}), divided into two sets: $\mathcal{F}_\mathrm{rand}=\lbrace{F_i}\rbrace_{i \in \mathcal{I}_\mathrm{rand}}$ and $\mathcal{F}_\mathrm{sing}=\lbrace{F_i}\rbrace_{i \in \mathcal{I}_\mathrm{sing}}$.  The 10 features in $\mathcal{F}_\mathrm{rand}$ are randomly built as in Test 2, with $\epsilon_i\sim \mathcal{U}(2e-3,1e-1)$, and $\bm{x}_i^C\sim  \mathcal{U}(-0.9,0.9)^2$, $i\in \mathcal{I}_\mathrm{rand}$. The position of the 9 features in $\mathcal{F}_\mathrm{sing}$ is instead given, in particular, for $\forall i\in \mathcal{I}_\mathrm{sing}$, $\bm{x}_j^C\in\lbrace-\frac{1}{2},0,\frac{1}{2}\rbrace\times\lbrace-\frac{1}{2},0,\frac{1}{2}\rbrace$, while $\epsilon_i\sim\mathcal{U}(6e-2,1e-1)$. For both set of features, i.e. for $i \in \mathcal{I}_\mathrm{rand}\cup \mathcal{I}_\mathrm{sing}$, $\theta_i\sim \mathcal{U}(0,2\pi)$, while the number of edges $n_i^e$ is chosen from a discrete uniform distribution between 4 and 16, as in Test 2. The specific feature data that were chosen for the proposed numerical experiment are reported in the Appendix.

On $\Omega=\Omega_0\setminus \lbrace \overline F,~F \in \mathcal{F}\rbrace$ we consider the following problem:
\begin{equation}
\begin{cases}
-\nabla \cdot(\kappa\nabla u)=0 &\text{in } \Omega\\
u=g_{\mathrm{D}} &\text{on } \Gamma_{\mathrm{D}}\\
\kappa\nabla u \cdot \bm{n}=0&\text{on } \Gamma_{\mathrm{N}}
\end{cases}
\end{equation}
with $\kappa$ being a discontinuous piecewise constant coefficient taking values $1$ and $100$ in a $4\times4$ chessboard pattern, as reported in Figure \ref{fig:diffusiongrid}. The features in $\mathcal{F}_{\mathrm{sing}}$ are hence exactly centered in the singularities produced by this piecewise constant coefficient. 
 The Dirichlet datum is chosen as a piecewise linear function. In particular, introducing
$$\mu_1(t)=\begin{cases}
2 &\text{if } t\in[-1,-\frac{1}{2})\\
-t+\frac{3}{2} &\text{if } t\in[-\frac{1}{2},0)\\
\frac{3}{2}&\text{if } t\in[0,\frac{1}{2})\\
-t+2&\text{if } t\in[\frac{1}{2},1],
\end{cases}\quad 
\mu_2(t)=\begin{cases}
-t &\text{if } t\in[-1,-\frac{1}{2})\\
\frac{1}{2}&\text{if } t\in[-\frac{1}{2},0)\\
-t+\frac{1}{2}&\text{if }t\in[0,\frac{1}{2})\\
0 &\text{if } t\in[\frac{1}{2},1],
\end{cases} $$
we define $$g_\mathrm{D}(x,y)=\begin{cases}
\mu_2(x) &\text{if }y=-1\\
\mu_1(-y)&\text{if }x=-1\\
\mu_1(x) &\text{if }y=1\\
\mu_2(-y)&\text{if }x=1.
\end{cases}$$
The initial mesh $\mathcal{T}_h^0$, which conforms the jumps in the coefficient $\kappa$, is reported in Figure \ref{fig:t3:initialMesh}. It is obtained by setting $h=1.41\cdot 10^{-1}$, ending up with $N=361$ degrees of freedom for variable $\uhs$

\begin{rem}
		Having introduced a coefficient $\kappa \neq 1$, we need to slightly modify the patch-local problems which we solve to obtain our (weakly) equilibrated flux reconstruction. For every $\bm{a}\in \N{h}$ we still solve a problem in the form of \eqref{eq:prob_eq_flux_nit}, with
	\begin{align*}
	&{m_{\bm{a}}(\bm{\rho}_h,\bm{v}_h):=(\kappa^{-1}\bm{\rho}_h,\bm{v}_h)_{\was}+\frac{1}{h_{\bm{a}}}\langle\kappa^{-1}\bm{\rho}_h \cdot \bm{n},\bm{v}_h \cdot \bm{n}\rangle_{\gamma_{\bm{a}}^\star}} \\
	&{L_{\bm{a}}(\bm{v_h}):=-(\psi_{\bm{a}}\nabla u_h^\star,\bm{v}_h)_{\was }-\frac{1}{h_{\bm{a}}}\langle \kappa^{-1}\psi_{\bm{a}}g,\bm{v}_h \cdot \bm{n}\rangle_{\gamma_{\bm{a}}^\star}}\\
	& R_{\bm{a}}(q_h)=(\psi_{\bm{a}}f-\nabla \psi_{\bm{a}}\cdot (\kappa\nabla u_h^\star),q_h)_{\was}+\langle \psi_{\bm{a}}g,q_h\rangle_{\gamma_{\bm{a}}^\star},
	\end{align*}
whereas the bilinear form $b_{\bm{a}}$ remains unchanged \cite{Smears_Vohralik_2020}.
We also need to change the definition of the numerical component of the estimator (see e.g. \cite{Vohralik2010}), by setting 
 $$(\estim{\sigma}^K)^2=||\kappa^{-\frac{1}{2}}(\flux +\kappa\nabla \uhs)||_{K\cap\Omegas}^2.$$ The use of this weighted norm derives from the need of bounding the energy norm of the overall error, which is defined as $||\kappa^\frac{1}{2}\nabla (u-\uhs)||_\Omega$.
\end{rem}

\begin{figure}
		\centering
			\includegraphics[width=0.45\linewidth]{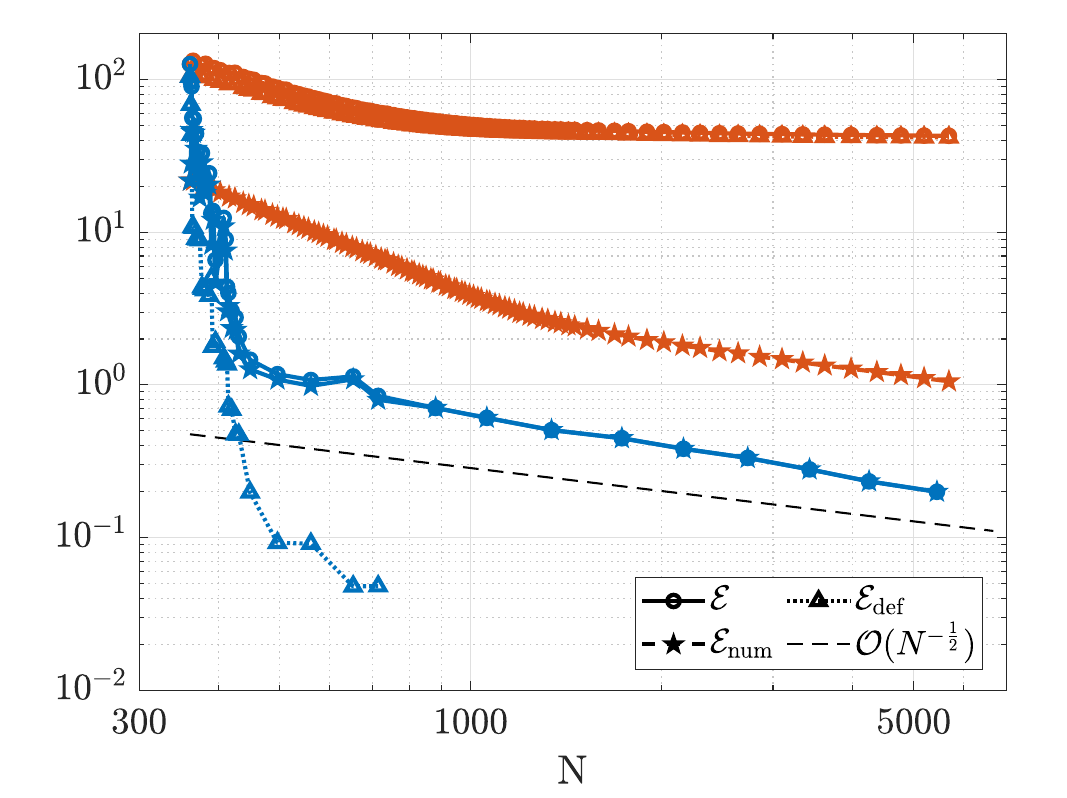}%
		\caption{Test 3: Convergence of the total estimator and of its components at the increase of the number of degrees of freedom. In red, no feature included; in blue, features included with $\alpha_3=1$.}
		\label{fig:t3:convergence}
\end{figure}

Figure \ref{fig:t3:convergence} reports the trend of the total error estimator and of its components during the adaptive process, considering both standard $h$-adaptivity and combined adaptivity. In both cases we choose $\theta=0.5$ in the D\"orfler marking strategy \eqref{Dorfler} and we stop before 6000 degrees of freedom are reached. The parameters $\alpha_1,\alpha_2,\alpha_3$ are all set equal to 1.
\begin{figure}
	\centering
	\begin{subfigure}{.5\textwidth}
		\centering
		\includegraphics[width=0.8\linewidth]{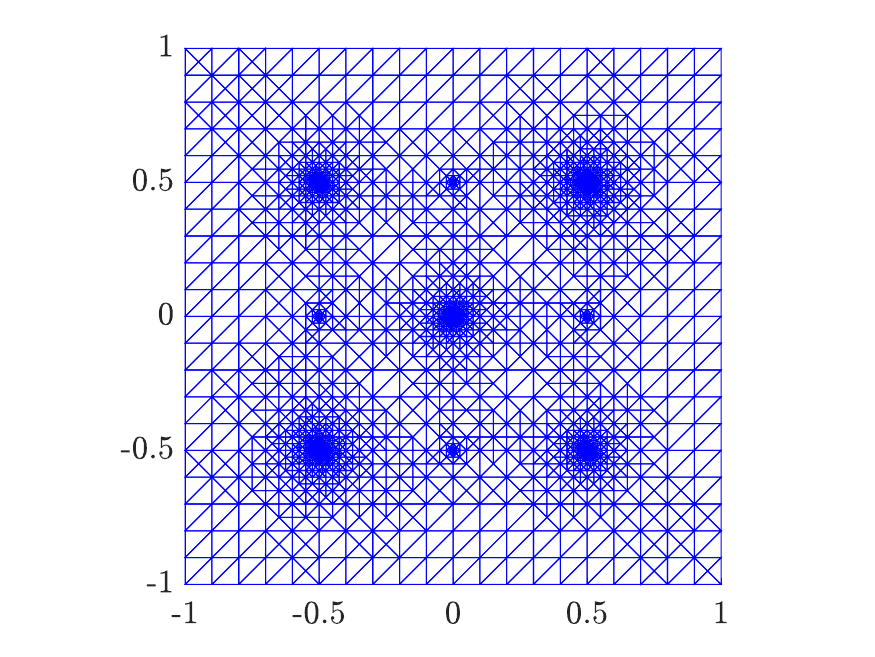}%
		\caption{Final mesh}
		\label{fig:t3:meshnoadd}
	\end{subfigure}%
	\begin{subfigure}{.5\textwidth}
		\centering
		\includegraphics[width=0.8\linewidth]{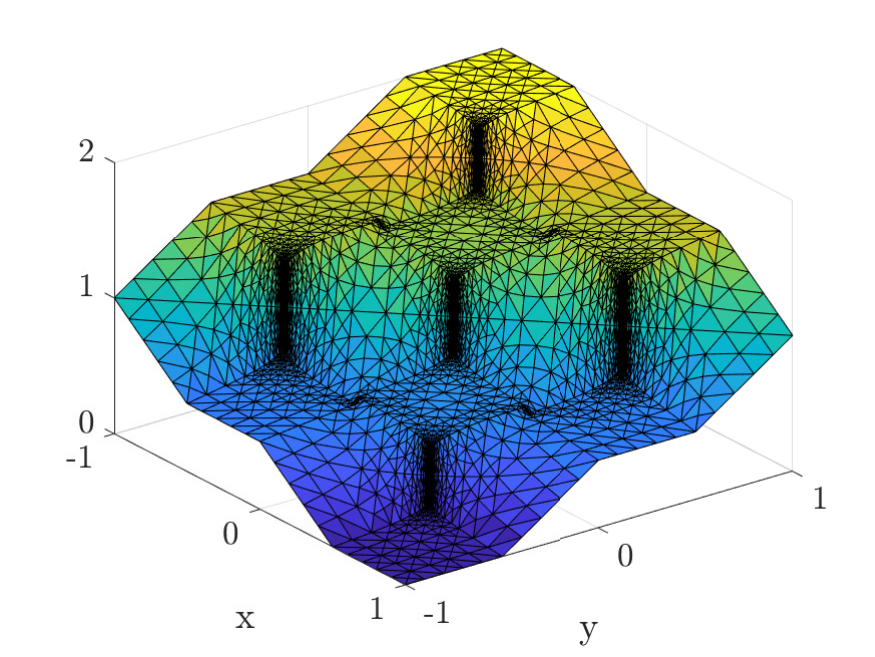}
		\caption{Final numerical solution}
		\label{fig:t3:meshnoadd_sol}
	\end{subfigure}
	\caption{Test 3: Final mesh and numerical solution in the case of standard mesh adaptivity (no feature is included).}
	\label{fig:t3:noadd}
\end{figure}
\begin{figure}[h]
	\centering
	\begin{subfigure}{.5\textwidth}
		\centering
		\includegraphics[width=0.8\linewidth]{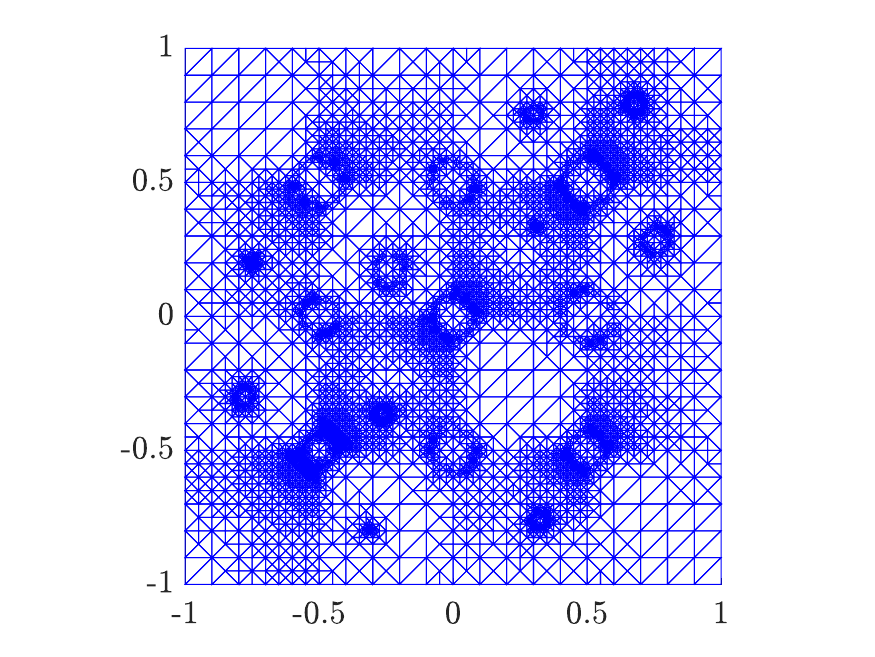}%
		\caption{Final mesh}
		\label{fig:t3:meshadd}
	\end{subfigure}%
	\begin{subfigure}{.5\textwidth}
		\centering
		\includegraphics[width=0.8\linewidth]{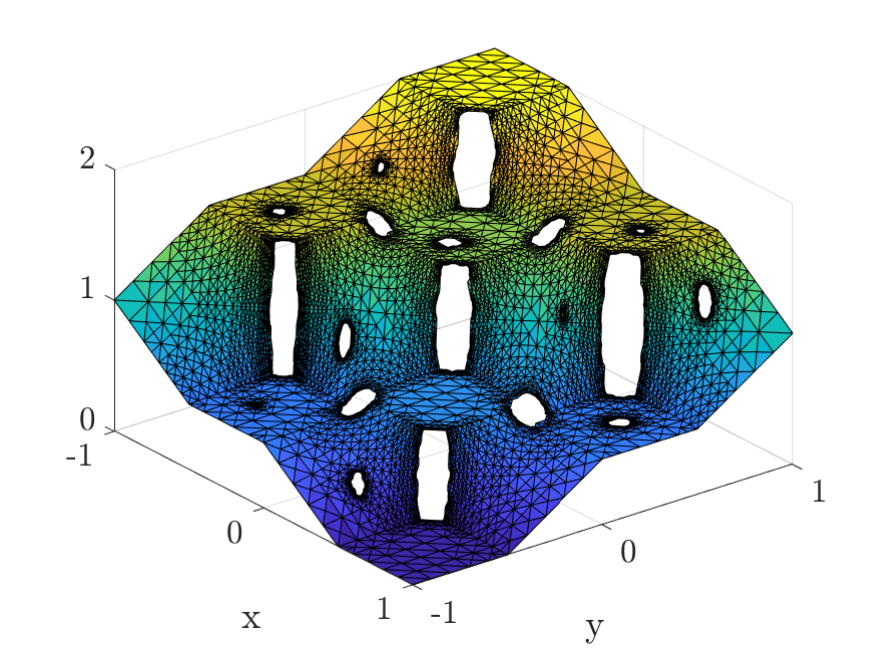}
		\caption{Final numerical solution}
		\label{fig:t3:meshadd_sol}
	\end{subfigure}
	\caption{Test 3: Final mesh and numerical solution in the case of combined adaptivity.}
	\label{fig:t3:meshes}
\end{figure}

 In the case of $h$-adaptivity we can observe how the total estimator converges very slowly, steered by the defeaturing component, which eventually reaches a plateau. The decrease of $\estim{\mathrm{def}}$ in the first iterations is related to the fact that the solution on the features centered in the singularities undergoes quick changes during the first refinement steps. The final mesh is reported in Figure \ref{fig:t3:meshnoadd}, while the final numerical solution is shown in Figure \ref{fig:t3:meshnoadd_sol}. In the case of the combined adaptivity we can instead observe how the total estimator tends to converge at the expected rate, completely steered by the numerical component, to which it is almost perfectly overlapped. The sudden decrease of the estimator in the first part of the adaptation is related to the fact that all the features in $\mathcal{F}_\mathrm{sing}$ are marked for inclusion during the first iterations. After that, in a few steps in which the mesh is refined close to these features, the singularities are finally excluded from the active mesh. The features in $\mathcal{F}_\mathrm{rand}$ are instead added a bit later, but all within iteration 27. The final mesh, along with the corresponding numerical solution is reported in Figure \ref{fig:t3:meshadd}.  

\begin{figure}
	\centering
	\includegraphics[width=0.45\linewidth]{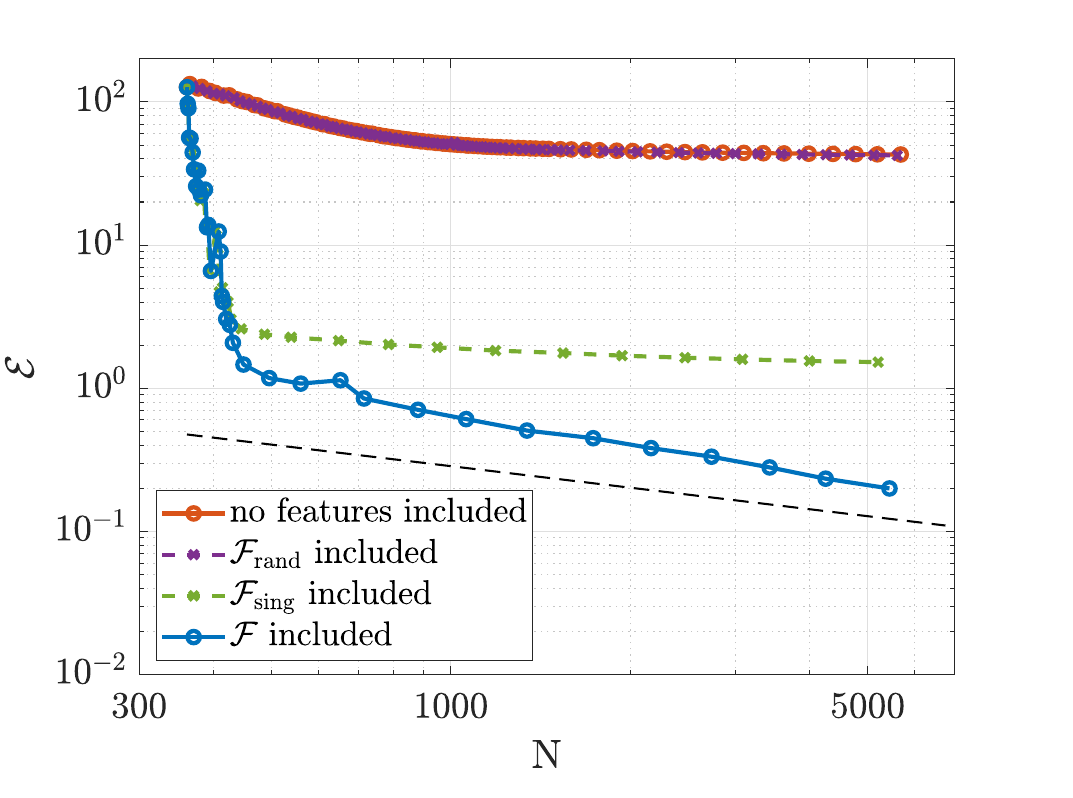}
	\caption{Test 3: Convergence with respect to the number of degrees of freedom (N) of the total estimator for different sets of features considered for inclusion.}
	\label{fig:t3:compare}
\end{figure}

Let us remark that, even if they are included later, the features in $\mathcal{F}_\mathrm{rand}$ are actually relevant for the accuracy of the solution. This is shown in Figure \ref{fig:t3:compare}, which reports the trend of the total estimator in a case in which the features in $\mathcal{F}_\mathrm{rand}$ are artificially not marked for inclusion although the estimator suggests to include them. Despite the rapid decrease related to the inclusion of the other features, $\estim{}$ tends to stagnate, not reaching the expected convergence rate. Figure \ref{fig:t3:compare} shows the trend of the estimator also in a case in which only the features in $\mathcal{F}_\mathrm{rand}$ are marked for inclusion, showing how, if the singularities are not excluded, there is almost no difference with the case in which the features are not added at all. This confirms the capability of the estimator, of course when the features are not artificially unmarked, to include the features in the most suitable order.

\section{Conclusions}
An adaptive strategy based on an equilibrated flux \emph{a posteriori} error estimator was proposed for defeaturing problems, extending the work in \cite{BCGVV} to the case of trimmed geometries. The focus is kept on Poisson's problem with Neumann boundary conditions on the feature boundary and, in particular, on the negative feature case. Similarly to \cite{BCV2022_arxiv} the adaptive strategy accounts both for standard mesh refinement and geometric refinement, i.e. for feature inclusion. The resulting adaptive strategy is called \emph{combined} adaptivity.

 In order to never remesh the considered domain, the features are progressively included into the geometry using {an unfitted finite element} strategy, and an equilibrated flux reconstruction is computed by imposing weakly, through Nitsche's method, the Neumann boundary conditions on the feature boundaries. Hence the estimator also accounts for the error committed in the imposition of such boundary conditions and for the consequent mass unbalance.

 The use of a reconstructed flux allows us to never evaluate the numerical flux on the boundary of the features, which is a great advantage for finite element discretizations, in which the gradient of the numerical solution is typically discontinuous across element faces. Furthermore, it {allows us to restrict} all contributions involving an unknown constant to cut elements and to the defeaturing component.
 
 The reliability of the estimator was proven in $\mathbb{R}^d$, $d=2,3$ and verified on three different numerical examples in $\mathbb{R}^2$, also involving multiple features of different sizes, both internal and on the boundary. 
 The generalization to the positive feature case is left to a forthcoming work, as well as the validation of the proposed adaptive strategy on three-dimensional geometries and on higher-order polynomial spaces, these two extensions having an impact mainly on computational aspects.

 \section*{Acknowledgments}
 Author Denise Grappein acknowledges that part of this work was carried out while holding a postdoctoral fellowship financed by the Italian Institute of Higher Mathematics (INdAM), hosted by the research unit of Politecnico di Torino. The work of Rafael Vázquez has been partially funded by Consellería de Educación, Ciencia, Universidades e Formación Profesional - Xunta de Galicia (2025-AD059, ED431F 2025/03 and ED431C 2025/09), and by the grant PID2024-156071NB-I00 funded by MICIU/AEI/10.13039/501100011033 and by ERDF/EU.
 
 The authors are grateful to Professor Claudio Canuto (DISMA, Politecnico di Torino) for some helpful discussions.
  \appendix
\section{Appendix 1}\label{appendix1}
Let us consider the asymmetric patch-local problem:
\emph{find $(\fluxa,\lam)\in \Mha \times \Qha$ such that}
\begin{equation}
\begin{cases}
m_{\bm{a}}(\fluxa,\bm{v}_h)-b_{\bm{a}}(\bm{v}_h,\lam)=L_{\bm{a}}(\bm{v}_h) & \forall \bm{v}_h \in \Mha\\
\hspace{2.2cm}b_{\bm{a}}^0(\bm{\fluxa},q_h)=R_{\bm{a}}^0(q_h)  & \forall q_h \in \Qha, \label{eq:prob_eq_flux_nit_asymm}
\end{cases}
\end{equation} where $m_{\bm{a}}$, $b_{\bm{a}}$ and $L_{\bm{a}}$ are defined in \eqref{eq:ma}, \eqref{eq:ba} and \eqref{Lh_nit} respectively, whereas $b_{\bm{a}}^0:\Mha\times\Qha \rightarrow \mathbb{R}$ and and $R_{\bm{a}}^0:\Qha \rightarrow \mathbb{R}$ are defined as
\begin{align*}
&b_{\bm{a}}^0(\bm{v}_h,q_h):=(q_h,\nabla \cdot \bm{v}_h)_{\was}\\
&R_{\bm{a}}^0(q_h):=(\psi_{\bm{a}}f-\nabla \psi_{\bm{a}}\cdot \nabla u_h^\star,q_h)_{\was}.
\end{align*}
For any $\bm{a} \in \mathscr{N}_h$ let $c_{\bm{a}}=-\overline{q_h}^\wa$, with $q_h \in Q_h(\wa)$. Since $q_h+c_{\bm{a}}\in \Qha$ we rewrite the second equation in \eqref{eq:prob_eq_flux_nit_asymm} as
\begin{equation}
(\nabla \cdot \fluxa,q_h+c_{\bm{a}} )_{\was }=(\psia f-\nabla \psi_{\bm{a}}\cdot \nabla u_h^\star,q_h+c_{\bm{a}})_{\was}.\label{eq:div_p0_asymm}
\end{equation}
According to \eqref{eq:hat_orth} and \eqref{eq:div_theo}, and remembering that Neumann boundary condition is imposed in strong form on the whole $ \partwa{\psi}\cap \Gamma_{\mathrm{N}}^0$
we can rewrite \eqref{eq:div_p0_asymm} as
\begin{align}
(\nabla \cdot \fluxa,q_h)_{\was}+\langle\fluxa \cdot \bm{n}+\psia g,c_{\bm{a}}\rangle_{\gamma_{\bm{a}}^\star}=(\psia f-\nabla \psi_{\bm{a}}\cdot \nabla u_h^\star,q_h)_{\was},\quad  \forall q_h \in Q_h(\wa) \label{eq:div_p0_intermediate}
\end{align}
which holds also for $q_h \in \mathcal{P}_\polydeg(K)$, since the polynomials in $Q_h(\wa)$ are discontinuous. We hence have that for any $K \in \mathcal{A}_h$ and for any $q_h \in \mathcal{P}_\polydeg(K) $
\begin{align}
(\nabla \cdot \flux,q_h)_{K\cap \Omegas}=\hspace{-0.3cm}
\sum_{\bm{a} \in \mathscr{N}_h(K)}(\nabla \cdot \fluxa,q_h)_{K\cap \Omegas}
&=(f,q_h)_{K \cap \Omegas}-\hspace{-0.3cm}\sum_{\bm{a} \in \mathscr{N}_h(K)}(\fluxa \cdot \bm{n}+\psia g,c_{\bm{a}})_{\gamma_{\bm{a}}^\star}\label{eq:div_p0_sum}
\end{align}
where we have exploited 
the fact that $\sum_{\bm{a} \in \mathscr{N}_h(K)}\nabla \psia=0$. 
{\section{Appendix 2: proof of Lemma 1}\label{appendix2}
	Following \cite[Appendix~A.1, Proof of Lemma~1]{GuzmanOlshanskii2018}, let
	$R_K : H^1(K) \rightarrow H^1(\mathbb{R}^d)$ be an extension operator such that
	$R_K v = v$ on $K$ and
	\begin{equation}\label{extension}
		\|R_K v\|_{\mathbb{R}^d} + h_K\|\nabla R_K v\|_{\mathbb{R}^d}
		\leq C_R\big(\|v\|_{K} + h_K\|\nabla v\|_{K}\big),\vspace{-0.3cm}
	\end{equation}
	where $C_R$ depends only on $d$ and on the shape regularity of the mesh, and not on $K$
	or $v$, and set $w:=R_Kv$.
	
	Let us first make explicit the dependence on $\delta_F$ of the trace inequality on $F$. Let us define a rescaled feature
	$\hat F := \delta_F^{-1}F$ having unit diameter, and for $\hat{\bm{x}} \in \hat F$ let
	$\hat w(\hat{\bm{x}}):=w(\delta_F\hat{\bm{x}})$. Since $\hat F$ is Lipschitz there exists a constant $C_{\hat F}>0$ such that {the following multiplicative trace inequality holds \cite[Theorem 1.6.6]{brenner_mathematical_2008}}:
	\begin{equation}\label{eq:trace_hat}
		\|\hat w\|^2_{\partial \hat F} \leq C_{\hat F}\big(\|\hat w\|^2_{\hat F}
		+ \|\hat w\|_{\hat F}\|\nabla \hat w\|_{\hat F}\big).
	\end{equation}
	Since $\hat{F}$ has unit diameter, the constant $C_{\hat F}$ depends only on $d$ and on the shape of $\hat F$.
	A change of variables gives
	\[
	\|w\|^2_{\partial F} = \delta_F^{d-1}\|\hat w\|^2_{\partial \hat F},\qquad
	\|w\|^2_{F} = \delta_F^{d}\|\hat w\|^2_{\hat F},\qquad
	\|\nabla w\|^2_{F} = \delta_F^{d-2}\|\nabla \hat w\|^2_{\hat F},
	\]
	so that, multiplying \eqref{eq:trace_hat} by $\delta_F^{d-1}$,
	\begin{equation}\label{eq:trace_F}
		\|w\|^2_{\partial F} \leq C_{\hat F}\big(\delta_F^{-1}\|w\|^2_{F}
		+ \|w\|_{F}\|\nabla w\|_{F}\big).
	\end{equation}
	Since $\gamma_K^\star \subseteq \partial F$ and $v = w$ on $K$, and using the
	arithmetic--geometric mean inequality
	$\|w\|_{F}\|\nabla w\|_{F} \leq \tfrac{1}{2}h_K^{-1}\|w\|^2_{F}
	+ \tfrac{1}{2}h_K\|\nabla w\|^2_{F}$ together with
	$\|\cdot\|_{F} \leq \|\cdot\|_{\mathbb{R}^d}$, we obtain from \eqref{eq:trace_F}
	\[
	\|v\|^2_{\gamma_K^\star} \leq C_{\hat F}\Big(\delta_F^{-1}\|w\|^2_{\mathbb{R}^d}
	+ \tfrac{1}{2}h_K^{-1}\|w\|^2_{\mathbb{R}^d}
	+ \tfrac{1}{2}h_K\|\nabla w\|^2_{\mathbb{R}^d}\Big).
	\]
	Since $h_K < h_0$ we have $\delta_F^{-1} < h_0\delta_F^{-1}h_K^{-1}$, whence
\begin{align*}
	\|v\|^2_{\gamma_K^\star} &\leq C_{\hat F}\big(\tfrac{1}{2}+h_0\delta_F^{-1}\big)
	\Big(h_K^{-1}\|w\|^2_{\mathbb{R}^d} + h_K\|\nabla w\|^2_{\mathbb{R}^d}\Big) \\ &{\leq
	 C_{\hat F}\big(\tfrac{1}{2}+h_0\delta_F^{-1}\big)h_K^{-1}
	\Big(\|w\|_{\mathbb{R}^d} + h_K\|\nabla w\|_{\mathbb{R}^d}\Big)^2}
\end{align*}
       {Now, using \eqref{extension} and the inequality}
	$(x+y)^2\leq 2(x^2+y^2)$ we obtain
	\[
	\|v\|^2_{\gamma_K^\star} \leq C_{\hat F}C_R^2\big(1+2h_0\delta_F^{-1}\big)
	\Big(h_K^{-1}\|v\|^2_{K} + h_K\|\nabla v\|^2_{K}\Big),
	\]
	which is \eqref{eq:loc_trace_lem}.
}
\section{Appendix 3}\label{appendice}
	In this appendix we provide the data of the features that were used for Test 2 and Test 3. A generic feature $F_i$ is assumed to be a polygon inscribed in a circle of radius $\epsilon_i$, centered in $\bm{x}_i^C=(x_i^C,y_i^C)$, having $n_i^e$ edges and rotated of an angle $\theta_i$ with respect to a reference orientation in which one of the vertexes is located in $\bm{x}_i^C+[0,\epsilon_i]$.
	\begin{table}[h]
	\renewcommand*{\arraystretch}{1}
	\centering 
	\caption{Feature data for Test 2\vspace{-0.3cm}}
	\label{test2_table}
	\begin{tabular}[t]{|c|ccccc||c|ccccc|}
		\hline
	$\bm{i}$ & $\epsilon_i$ & $x_i^C$ &$y_i^C$ & $n_i^e$ & $\theta_i$ & $\bm{i}$ & $\epsilon_i$ & $x_i^C$ &$y_i^C$ & $n_i^e$ & $\theta_i$ \\
	\hline
	1  & 0.0128 & 0.3695 & 0.6692 & 4  & 224.9859 & 20 & 0.0156 & 0.8396 & 0.8958 & 5  & 154.4009 \\
	2  & 0.0239 & 0.1714 & 0.7986 & 6  & 24.2550  & 21 & 0.0057 & 0.8600 & 0.4401 & 8  & 134.2678 \\
	3  & 0.0314 & 0.5347 & 0.7247 & 11 & 29.1858  & 22 & 0.0367 & 0.2399 & 0.4603 & 16 & 73.3359 \\
	4  & 0.0277 & 0.1893 & 0.1487 & 10 & 190.0490 & 23 & 0.0045 & 0.5668 & 0.4321 & 4  & 68.2234 \\
	5  & 0.0352 & 0.8916 & 0.2945 & 7  & 262.3478 & 24 & 0.0158 & 0.3045 & 0.7608 & 7  & 52.4692 \\
	6  & 0.0487 & 0.6498 & 0.1642 & 16 & 125.5876 & 25 & 0.0202 & 0.2724 & 0.8842 & 15 & 19.4049 \\
	7  & 0.0269 & 0.6300 & 0.5445 & 8  & 87.1778  & 26 & 0.0388 & 0.7788 & 0.5486 & 14 & 108.5075 \\
	8  & 0.0174 & 0.1094 & 0.3495 & 13 & 56.1573  & 27 & 0.0023 & 0.5995 & 0.2871 & 4  & 185.3529 \\
	9  & 0.0110 & 0.1255 & 0.6014 & 8  & 126.4299 & 28 & 0.0302 & 0.0000 & 0.6995 & 11 & 6.0078 \\
	10 & 0.0156 & 0.2704 & 0.6302 & 11 & 81.7039  & 29 & 0.0386 & 0.0000 & 0.3592 & 6  & 357.7583 \\
	11 & 0.0422 & 0.4453 & 0.8543 & 16 & 184.2681 & 30 & 0.0452 & 0.0000 & 0.8618 & 6  & 336.7216 \\
	12 & 0.0058 & 0.2489 & 0.3096 & 11 & 124.2501 & 31 & 0.0475 & 0.0000 & 0.1621 & 6  & 75.1260 \\
	13 & 0.0148 & 0.6500 & 0.8789 & 4  & 115.4929 & 32 & 0.0036 & 0.0000 & 0.5210 & 9  & 238.0898 \\
	14 & 0.0254 & 0.4708 & 0.3901 & 9  & 281.4964 & 33 & 0.0484 & 1.0000 & 0.8710 & 7  & 265.5165 \\
	15 & 0.0138 & 0.3684 & 0.3172 & 8  & 346.0244 & 34 & 0.0249 & 1.0000 & 0.2814 & 16 & 233.4964 \\
	16 & 0.0492 & 0.4653 & 0.1502 & 7  & 249.4054 & 35 & 0.0295 & 1.0000 & 0.1099 & 7  & 275.1374 \\
	17 & 0.0045 & 0.3452 & 0.1051 & 5  & 17.9789  & 36 & 0.0273 & 1.0000 & 0.5045 & 5  & 147.7730 \\
	18 & 0.0026 & 0.6548 & 0.6676 & 10 & 66.5531  & 37 & 0.0376 & 1.0000 & 0.7132 & 13 & 176.2045 \\
	19 & 0.0247 & 0.8788 & 0.6407 & 7  & 48.8767 & & & & & & \\
	\hline
	\end{tabular}%
\end{table}
	\begin{table}[h]
	\renewcommand*{\arraystretch}{1}
	\centering 
	\caption{Feature data for Test 3\vspace{-0.3cm}}
	\label{test3_table}
	\begin{tabular}[t]{|c|ccccc||c|ccccc|}
		\hline
		$\bm{i}$ & $\epsilon_i$ & $x_i^C$ &$y_i^C$ & $n_i^e$ & $\theta_i$ & $\bm{i}$ & $\epsilon_i$ & $x_i^C$ &$y_i^C$ & $n_i^e$ & $\theta_i$ \\
		\hline
	1  & 0.0617 & -0.5000 & -0.5000 & 8  & 87.1778 & 11 & 0.0219 & -0.7489 & 0.2029  & 4  & 115.4929 \\
	2  & 0.0674 & -0.5000 & 0.0000  & 13 & 56.1573  & 12 & 0.0297 & 0.3001  & 0.7578  & 10 & 70.8553 \\
	3  & 0.0915 & -0.5000 & 0.5000  & 8  & 126.4299 & 13 & 0.0276 & -0.2632 & -0.3656 & 16 & 137.5424 \\
	4  & 0.0830 & 0.0000  & -0.5000 & 14 & 26.4852  & 14 & 0.0090 & -0.3096 & -0.7899 & 7  & 85.2644 \\
	5  & 0.0791 & 0.0000  & 0.0000  & 11 & 338.4382 & 15 & 0.0051 & 0.3097  & 0.3352  & 4  & 165.5202 \\
	6  & 0.0833 & 0.0000  & 0.5000  & 11 & 81.7039  & 16 & 0.0495 & 0.7576  & 0.2813  & 13 & 327.8332 \\
	7  & 0.0719 & 0.5000  & -0.5000 & 14 & 261.8490 & 17 & 0.0313 & 0.6791  & 0.7915  & 12 & 39.8837 \\
	8  & 0.0981 & 0.5000  & 0.0000  & 14 & 212.4525 & 18 & 0.0315 & 0.3245  & -0.7584 & 6  & 243.2179 \\
	9  & 0.0921 & 0.5000  & 0.5000  & 16 & 184.2681 & 19 & 0.0585 & -0.2369 & 0.1675  & 12 & 225.4836 \\
	10 & 0.0347 & -0.7811 & -0.3010 & 11 & 124.2501 & & & & & & \\
	\hline
	\end{tabular}%
\end{table}

\printbibliography
\end{document}